\documentclass[preprint,11pt]{article}

\usepackage{braket,amsfonts}
\usepackage{amssymb}
\usepackage{amsthm}
\usepackage{array}
\usepackage{url}

\usepackage{pgfplots}
\usepackage{hyperref}
\usepackage{authblk}
\usepackage[a4paper, total={7in, 9in}]{geometry}
\usepackage{graphicx}

\usepackage{amsopn}

\tikzset{
  treenode/.style = {shape=rectangle, rounded corners,
                     draw, align=center,
                     top color=white, bottom color=blue!20},
  root/.style     = {treenode, font=\Large, bottom color=red!30},
  env/.style      = {treenode, font=\ttfamily\normalsize},
  dummy/.style    = {circle,draw}
}
\usetikzlibrary{positioning,calc}

\usepackage{mathtools}

\newtheorem{remark}{Remark}[section]
\newcommand{\R}{\mathbb{R}}

\usepackage{xspace}
\usepackage{bold-extra}
\usepackage[most]{tcolorbox}

\usepackage{subcaption}
\usepackage{wrapfig}

\renewcommand{\R}{\mathbb{R}}

\newcommand{\cJ}{\mathcal{J}}

\newcommand{\T}{{\mathbb T}}

\usepackage[ruled,vlined]{algorithm2e}
\SetKwComment{Comment}{}{} 

\pgfplotsset{compat=1.16}
\begin{document}


\title{Computation and Control of Unstable Steady States for Mean Field  Multiagent Systems}

\author{Sara Bicego}
\author{Dante Kalise}
\author{Grigorios A. Pavliotis}
\affil{Department of Mathematics, Imperial College London, United Kingdom. 
		\texttt{e-mail:\{s.bicego21,dkaliseb,g.pavliotis\}@ic.ac.uk}}

\maketitle
\begin{abstract}
We study McKean-Vlasov PDEs obained as the mean field limit of interacting particle systems driven by noise, modeling phenomena such as opinion dynamics. We are interested in systems that exhibit phase transitions i.e. non-uniqueness of stationary states for the corresponding PDE, in the mean field limit. We develop an efficient numerical scheme for identifying all steady states (both stable and unstable) of the mean field McKean-Vlasov PDE, based on a spectral Galerkin approximation combined with a deflated Newton's method to handle the multiplicity of solutions. Having found all possible equilibra, we formulate an optimal control strategy for steering the dynamics towards a chosen unstable steady state. The control is computed using iterated open-loop solvers in a receding horizon fashion. We demonstrate the effectiveness of the proposed steady state computation and stabilization methodology on several examples, including the noisy Hegselmann-Krause model for opinion dynamics and the Haken-Kelso-Bunz model from biophysics. The numerical experiments validate the ability of the approach to capture the rich self-organization landscape of these systems and to stabilize unstable configurations of interest. The proposed computational framework opens up new possibilities for understanding and controlling the collective behavior of noise-driven interacting particle systems, with potential applications in various fields such as social dynamics, biological synchronization, and collective behavior in physical and social systems.

\end{abstract}

\tableofcontents

\section{Introduction} Multiagent systems, widely used as models in biology, epidemics, ecology, and the social sciences, exhibit collective behavior and coherent structures at the macroscale due to interactions at the microscale \cite{IPS_crowd,IPS_epidemic,IPS_social, toscani2014,Motsch_Tadmor}. Examples include swarming, flocking, bacteria synchronization, socioeconomic and life sciences~\cite{Naldi-al-2010, Degond_al_2014,ABCB} and clustering emergence in opinion dynamics~\cite{helfmann2023modelling}. In this paper we focus on multiagent systems that are described in terms of a system of weakly interacting diffusion processes. The analysis of such stochastic interacting agent systems is facilitated by passing to the mean field limit, where the $N-$particle Fokker-Planck equation is replaced by an equation for the evolution of the $1-$particle distribution function, the nonlinear and nonlocal McKean-Vlasov PDE~\cite{FrankFP, sznitman1991}. This paradigm has proven to be very effective in a variety of applications, e.g. physics, astrophysics and biology \cite{FP_Black_Holes, FP_review}, modelling of traffic flow and opinion formation  \cite{FP_traffic,FP_opinions}, and convergence analysis of Artificial Neural Networks \cite{FP_NNs,IPS_NN1,IPS_NN2,IPS_NN3,IPS_NN4}, to name but a few. Emergent patterns in these systems are represented by stationary states of the McKean-Vlasov PDE, where a common feature is the coexistence of various stable and unstable configurations and transitions between them. The nature and number of steady states depends on the interaction and noise strengths, measuring the interplay between interaction and diffusion~\cite{Dawson1983, shiino1987, CGPS2020}.The transition between different macroscopic states can be interpreted as a disorder/order phase transition. The class of multiagent models that we consider in this paper are characterized by the presence of multiple stationary states~\cite{Degaldino, FrankFP,BCCD}.  The first main objective of this paper is to develop accurate and efficient computational methodologies for identifying all stationary states, both stable and unstable, of noise-driven interacting particle systems in their mean field limit.

The multiagent systems we study possess a very rich self-organization landscape. Stable equilibria often exhibit large basins of attraction, however, convergence may be slow depending on confinement and interaction potentials~\cite{monmarche2024local}. In the context of linear Fokker-Planck type PDEs with a unique Gibbsian stationary state, it is of interest to introduce external forcing terms which can accelerate convergence towards globally stable steady states~\cite{Kunisch_al_2018, Lelievre_al_2013}.  This is naturally addressed in the framework of nonlinear control theory and feedback stabilization. In this paper we follow this control-theoretical approach, and, as a second main objective, we study the optimal stabilization of McKean-Vlasov PDEs towards unstable steady states.  As a prototypical example, we will focus on opinion dynamics such as the noisy Hegselmann-Krause model, where the goal is to promote or prevent the emergence of polarization or consensus in the population, and where these steady states undergo a phase transition~\cite{Chazelle_al2017a, GPY2017, CGPS2020,during_wolfram24}. Optimal control methodologies have been applied to the problem of synchronization for Kuramoto-type models \cite{control_Kuramoto}. The methodology developed in this paper applies in a straightforward manner to this problem; in fact we study a variant of the noisy Kuramoto model, with an external potential, in Sec~\ref{subsec:O2}. Applying our methodology in order to achieve synchronization for the noisy Kuramoto model is an easier problem to the Heglselmann-Krause one, since the Kuramoto model is characterized by a continuous phase transition, as opposed to the discontinuous phase transition that is present in bounded confidence models.

\paragraph{Models, methods, and related literature.} We study multiagent systems with weakly interacting diffusions. For such systems the passage to the mean field limit~\cite{sznitman1991, oelschlager1984} and the structure and nature of phase transitions for the mean field dynamics is well understood~\cite{Dawson1983, ChayesPanferov2010, CGPS2020}. We consider McKean-Vlasov PDEs characterized by confinement and interaction potentials. Similarly to~\cite{SGGPUV2019}, see also~\cite{pearson22, Chavanis2014}, we discretize the stationary McKean-Vlasov PDE using a spectral Galerkin method,  approximating the steady state solution as a truncated Fourier series. This translates into a system of nonlinear equations for the Fourier coefficients which can be written in vector form. The integrals involved in the discretized operators are computed using Gaussian quadrature. The constraint that mass is conserved, i.e. that the solution of the McKean-Vlasov PDE integrates to 1, is imposed by adding  a constraint to the system. Positivity of the density is enforced by filtering out solutions with negative values. 

{\color{black} Finding multiple solutions to nonlinear PDEs such as the Allen-Cahn or the Yamabe equation has been previously studied in \cite{deflation}, where an abstract deflated Newton method in Banach spaces was proposed. The same method has been applied for bifurcation analysis in \cite{deflation_2}.} 
The standard Newton's method is combined with a deflation technique to systematically eliminate known solutions from consideration. The PDE residual is iteratively modified so that the roots of the deflated operator coincide with the remaining unknown solutions. The deflation operator is constructed using a shifted exponential-norm deflation operator. The derivative of the deflated residual, required for Newton's method, is efficiently computed using a recursive formula. The deflation loop continues until the root-finding procedure diverges for a given deflated operator and initial guess. This algorithm yields a complete taxonomy of the equilibrium points of the system, for which we formulate an optimal stabilization problem.

Optimal control for interacting multiagent systems has been extensively studied over the last decade across different scales. At the microscopic level, the work \cite{cfpt} deals with consensus control over the Cucker-Smale model with $N$ microscopic agents. This idea has been extended to the kinetic level in \cite{albi2014kinetic,evac} where the control acts over a Boltzmann type PDE. In this paper, we focus on the control of densities at the mean field scale, that we formulate as a PDE-constrained optimization problem for the McKean-Vlasov PDE. First-order optimality conditions for this problem have been previously studied in \cite{PMP_PDE_Aronna,robfleig} for the Fokker-Planck equation, and in \cite{acfk,mfpmp,pearson22} for the McKean-Vlasov PDE. Although related, our work is conceptually different as it focuses on the mean-field/McKean-Vlasov nonlinear and nonlocal PDEs where multiple steady states can exist, and which cannot be stabilized using open-loop controllers. To circumvent this difficulty, we resort to the synthesis of feedback laws \cite{Kunisch_al_2018,dkk_21} via nonlinear model predictive control (MPC) \cite{borzisurvey,mfgmpc,mdmpc,AKMPC}. We also mention the (formal) link between optimal control methodologies for mean field dynamics and linear response theory, e.g. the study of the response of the system to a weak external forcing (control), in particular close to the transition point~\cite{LPZ2020, LPZ2021}. It was shown in these references that linear response theory breaks down close to the transition point, and that a weak external forcing can have a substantial impact on the behaviour of the mean field dynamics.  
{\color{black} We emphasize the fact that the goal of this paper is to propose and test a systematic and robust methodology for identifying {\bf all} (stable and unstable) steady states of the McKean-Vlasov PDE, and of developing optimal control methodologies for steering the dynamics towards a chosen (unstable) steady state. The rigorous analysis of the methodology developed in this paper is a topic for future work.}

In summary, the main contributions of this paper are:
\begin{enumerate}
\item The development of a numerical scheme for identifying all steady states (up to translation invariance), both stable and unstable, of noise-driven interacting particle systems in their mean field limit, posed on the torus. The proposed method combines a spectral Galerkin approximation of the McKean-Vlasov PDE with a deflated Newton's method to handle the multiplicity of solutions. This allows for a comprehensive characterization of the self-organization landscape of these systems, including the coexistence of various stationary configurations and the transitions between them.

\item The formulation and solution of an optimal control problem for stabilizing the McKean-Vlasov PDE towards an unstable steady state. The control enters as an external forcing term in the PDE and is computed by solving a finite-horizon optimization problem using first-order optimality conditions. The control is synthesized using MPC, enabling the stabilization of unstable stationary solutions in a receding horizon fashion.

\item The application of the proposed methodology to several examples of interacting particle systems exhibiting phase transitions and multiple steady states, including the noisy Hegselmann-Krause model for opinion dynamics and the Haken-Kelso-Bunz model from biophysics. The numerical results demonstrate the effectiveness of the approach in capturing the rich self-organization behavior of these systems and stabilizing unstable configurations of interest. Detailed numerical studies illustrate the role of noise and interaction strength in shaping the landscape of stationary solutions in the considered models. We analyze the bifurcations and phase transitions in the number and stability of steady states as the system parameters are varied, providing insights into the collective behavior of these systems.
\end{enumerate}

\paragraph{Outline of this paper.}The rest of the paper is organized as follows. In Section~\ref{sec:models} we introduce the class of models under consideration, from the Brownian particles dynamics to the Fokker-Planck equation arising at the mean field limit. Section~\ref{sec:comp_multiple_states} covers the identification of stationary states at the macroscopic level. We fix the spatial discretization of the steady state PDE using a spectral Fourier Galerkin method. The resulting system of nonlinear equations is solved via Newton's method combined with a deflation algorithm to handle the multiplicity of solutions typical of this class of problems. We also present two numerical examples in which our proposed steady-states-finding-procedure is tested and validated against well-known linear stability analysis and bifurcation results. In Section $4$, we modify the evolutionary Fokker-Planck equation to include an external control force. This control signal is designed to stabilize the system at a desired unstable stationary state through the optimization of a suitable cost functional. The resulting finite-horizon optimal control problem is addressed at the discrete level using first-order optimality conditions, and the control action is synthesized in feedback form via nonlinear model predictive control. In Section~\ref{sec:numerics}, we test this stabilization routine for two different models across one and two-dimensional domains. In Section~\ref{sec:conclusions} we present conclusions and outline future work. {\color{black}Finally, we provide an appendix in Section~\ref{sec:appendix}, where we address some more theoretical points regarding our numerical methodology.}

\section{The McKean-Vlasov PDE and Characterization of Steady States}
\label{sec:models}

We consider a system of $N$ interacting particles driven by additive noise:
\begin{equation}\label{SDE}
    d X^i_t = - \nabla V(X^i_t) \,dt - \dfrac{1}{N} \sum\limits_{j=1}^N \nabla W (X^i_t,X^j_t) \,dt + \sqrt{2 \beta^{-1}} \,dB^i_t\,, \quad X_0^i \sim \nu \; ,\; i=1,\dots,N,
\end{equation}
where $X^i_t$ denotes the state of the $i-$th agent at time $t$. Here, $V$ and $W$ denote the confining and interaction potential, respectively, $\beta$ denotes the inverse temperature, and $B_t:=\big(B^1_t,\dots,B^N_t\big)$ is a standard $N-$dimensional Brownian motion. The confining potential is a symmetric function, $W(x,y) = W(y,x)$.  We consider the dynamics on the torus $\T^d$, in one and two dimensions. The noisy Hegselmann-Krause model~\cite{Chazelle_al2017a, GPY2012, GGSP2021} that we will study later on is a prototypical example of a system of the form~\eqref{SDE}. The white noise in~\eqref{SDE} models uncertainty and it can be easily replaced by noise processes with a non-trivial spatio-temporal structure~\cite{SGGPUV2019}.
As is well known~\cite[Ch. 4]{Pavl2014} that the dynamics~\eqref{SDE} has a gradient structure with respect to the Hamiltonian  
\begin{equation}\label{e:hamiltonian}
    \mathcal{H}^N(X^1,\dots,X^N) := \sum\limits_{i=1}^N V(x^i) + \dfrac{1}{2 N}\sum\limits_{i,j=1}^N W(x^i,x^j)
\end{equation}
and reads
\begin{equation}\label{e:IPS}
    d X^i_t = - \nabla \mathcal{H}^N(X^1_t,\dots,X^N_t) \,dt + \sqrt{2 \beta^{-1}}\,dB^i_t\,.
\end{equation}
The $N-$particle distribution function   $\rho^N(t,x_1, x_2, \dots , x_N)$ satisfies the Fokker-Planck equation associated with the $N-$particle system as
\begin{equation*}\label{n_part}
    \partial_t \rho^N_t =  \nabla_N\cdot\bigg(\nabla_N\mathcal{H}^N\rho^N_t\bigg) + \beta^{-1} \Delta_N\rho^N_t\,.
\end{equation*}
Under appropriate assumptions on the Hamiltonian and certainly when considering the dynamics on the torus with smooth confining and interaction potentials,  the $N-$particle Fokker-Planck equation has a unique steady state given by the Gibbs measure 
\begin{equation}\label{e:gibbs_N}
    \rho^N_\infty := \dfrac{1}{Z^N_\beta} e^{-\beta\mathcal{H}^N(x_1,\dots,x_N)},\qquad \quad Z^N_\beta := \int e^{-\beta\mathcal{H}^N(x_1,\dots,x_N)} \,  d x_1 \dots x_N \,.
\end{equation}
    We refer to~\cite[Ch. 4]{Pavl2014} for details.

For exchangeable interacting particle systems and with chaotic initial data, and when the interaction strength is inversely proportional to the number of particles as in~\eqref{e:hamiltonian}, we can pass to the mean field limit $N \rightarrow +\infty$. In particular, the $N-$particle distribution function becomes approximately the product of $1-$particle distribution functions, $\rho^{N}(x_1,\dots, x_N,t) \approx \Pi_{\ell =1}^N \rho(x_{\ell},t)$; the $1-$particle distribution function satisfies the {\emph McKean-Vlasov} PDE
\begin{equation}\label{fp1}
    \partial_t \rho = \nabla \cdot\bigg(\rho\big[(\nabla W*\rho) + \nabla V \big]    \bigg) + \beta^{-1} \Delta\rho\,,
\end{equation}
in the domain $\T^d$ with periodic boundary conditions and initial distribution $\rho(0,x) = \rho_0$, with
\begin{equation*}\label{e:star}
    \big(W*\rho\big)\,(x) := \int\limits_{\T^d} W(x,y)\rho(y)dy\,.
\end{equation*}
We refer to the classic references~\cite{sznitman1991, oelschlager1984} for proofs and refinements of the above and to~\cite{DGPS23, Lacker_LeFlem_2023} for modern approaches to entropic propagation of chaos.

Unlike the $N-$particle system~\eqref{e:IPS} that is always reversible with respect to the unique Gibbs invariant measure~\eqref{e:gibbs_N}, it is well known that the mean field dynamics can exhibit phase transitions, i.e. non-uniqueness of stationary states~\cite{Dawson1983, shiino1987, FrankFP, CGPS2020}. This is the case, for example, for the mean field dynamics in unbounded domains and for a non-convex confining potential and a quadratic interaction potential at sufficiently low temperatures~\cite{Dawson1983, Tugaut2014}. For the dynamics on the torus, which is the case that we will consider in the paper, the concept of convexity is replaced by that of \emph{$H-$Stability}~\cite{ruelle1999statistical, ChayesPanferov2010, CGPS2020}: The mean field dynamics exhibits a unique stationary state at all temperatures if and only if the interaction potential has no non-negative Fourier coefficients. When the interaction potential has negative Fourier coefficients, then a phase transition occurs at sufficiently low temperatures.  We refer to~\cite{ChayesPanferov2010, CGPS2020} for a detailed analysis of phase transitions in the absence of a confining potential and to~\cite{Tamura1984, Degaldino, bertoli2024} for partial results when a confining potential is present on the torus. It is important to note that, in the absence of a confining potential, the uniform distribution on the torus is always a stationary state. Furthermore, past the phase transition, steady states are translation invariant~\cite{Giacomin_al_2010}. We can identify a particular steady state by choosing \color{black}appropriately the initial condition~\cite{bertoli2024}. \color{black} Translation invariance is removed by the addition of a confining potential, e.g. a "magnetic field"~\cite{Degaldino}.  

Several equivalent characterizations of stable and unstable steady states of \eqref{fp1} exist: they can be defined as solutions to the steady state PDE
\begin{equation}\label{stationary_FP}
    \beta^{-1} \Delta\rho+\nabla \cdot\bigg(\rho\big[(\nabla W*\rho)  + \nabla V\big]\bigg) = 0\,,
\end{equation}
as critical points of the free energy functional 
\begin{equation}\label{freeE}
    \mathcal{F}[\rho] = \beta^{-1}\int_{\T^d} \rho \log(\rho)\,dx + \int_{\T^d} V \rho\, dx +\frac{1}{2}\int_{\T^d} \int_{\T^d} W(x,y) \rho(x) \rho(y) \, dx dy\,,
\end{equation}
or fixed points of the Euler-Lagrange equations for the free energy functional~\cite{Bavaud_1991} 
\begin{equation}\label{e:integ_eqn}
\rho = \frac{1}{Z} e^{-\beta (V + W \ast \rho)}, \quad Z:= \int_{\T^d} e^{-\beta (V + W \ast \rho)} \, dx.
\end{equation}
Alternatively, it can be shown that this integral equation, the \emph{Kirkwood-Monroe map} is equivalent to the stationary McKean-Vlasov PDE~\eqref{fp1}; see~\cite[Thm 4.1]{Tamura1984} or~\cite[Thm. 1]{Dressler1987}. {\color{black} A rigorous proof of the equivalence between steady states obtained through the previously mentioned methods is presented in~\cite[Prop 2.4]{CGPS2020}. Furthermore, ~\cite[Thm 2.3]{CGPS2020} establishes the regularity of steady state solutions as $\rho \in H^1 \cap \mathcal{P}_{ac}(\Omega)$, where $\mathcal{P}_{ac}(\Omega)$ denotes the space of absolutely continuous probability measures over $\Omega$. The solutions are in fact strictly positive \cite{nickl2024bayesiannonparametricinferencemckeanvlasov}.}

We reiterate that finding critical points of the free energy functional~\eqref{freeE}, roots of the nonlinear map $F_{\kappa}[\rho] := \rho - \frac{1}{Z} e^{-\beta (V + W \ast \rho)}$ in~\eqref{e:integ_eqn} and stationary solutions of the McKean-Vlasov PDE are \emph{equivalent}~[Prop 2.4]\cite{CGPS2020}.\footnote{In fact, these statements are also equivalent to finding global minimizers of the entropy dissipation functional, but this observation will not be used in this paper.} One of the main goals of this paper will be the development of efficient numerical methodologies that enable us to calculate all solutions of the Kirkwood-Monroe equation, i.e. all critical points of the free energy functional.

\section{Computing multiple steady states in the McKean-Vlasov PDE}
\label{sec:comp_multiple_states}

The variety of applications of this problem has motivated a flourishing literature on numerical solvers for \eqref{fp1}. For a review of finite difference and finite element methods, we refer the reader to \cite{fd_vs_fe}, beyond which we cite the structure preserving scheme proposed in \cite{pareschi18}, where the positivity of $\rho$ and its mass preservation are guaranteed by construction. Alternative solvers include Spectral Methods \cite{spectral_FP,spectral_FP2}, Deep Learning \cite{DeepFP,FP_PINNs,neural_FP,DeepFP2}, and Tensor Train approximation \cite{TT_FP}.

In this paper we are concerned with the identification of invariant probability measures of \eqref{SDE} obtained as solutions of the stationary McKean-Vlasov PDE~\eqref{stationary_FP}. The proposed numerical procedure relies on a spectral approximation of $\rho$ in terms of $L$ Fourier modes, for which \eqref{fp1} can be written as a bilinear system of $L$ equations. In these settings, the identification of steady states reduces to a root finding problem, addressed via Newton's method.  In order to identify distinct solutions, typically one applies the routine multiple times starting from many different initial guesses, and relying on the assumption that those starting points lie in different basins of attraction. Here instead, we rely on the infinite-dimensional deflation algorithm proposed in \cite{deflation}, where the PDE residual is iteratively modified in such a way to eliminate known solutions from consideration. 
Alternative approaches include numerical continuation \cite{continuation,continuation2}, and numerical integration of the Davidenko differential equation associated with the original nonlinear problem \cite{DDE}.
%
%
\subsection{Spectral Galerkin method for the Stationary Fokker--Planck equation}
We consider a Galerkin approximation of the stationary Fokker--Planck equation \eqref{stationary_FP}. Its weak form over a periodic domain $\Omega$ reads
\begin{equation*}
     \beta^{-1}\int\limits_{\Omega} \Delta\rho \,\varphi \, dx + \int\limits_{\Omega}\nabla \cdot\bigg(\rho\big[(\nabla W*\rho) +\nabla V\big]\bigg)\,\varphi \, dx\, = 0,
\end{equation*}
with test functions $\varphi\in \mathcal{V}:=\{\varphi:\Omega\to\mathbb{R}\;\big|\; \text{$\varphi$ smooth, periodic} \}$. After integrating by parts and under the periodicity assumption for both $\rho$ and $\varphi$, this leads to
\begin{equation*}
    0 = - \beta^{-1}\bigg( \int\limits_{\Omega} \nabla\rho\cdot\nabla\varphi\,dx\bigg) - \bigg(\int\limits_{\Omega} \big(\rho\big[(\nabla W*\rho) + \nabla V\big]\cdot\nabla\varphi dx\bigg)\,,\qquad \forall \varphi \in \mathcal{V}\,.
\end{equation*}
We consider a global spectral basis $\{\psi_i\}_{i=1}^L$  spanning an $L-$dimensional subspace $\mathcal{V}_L\subset\mathcal{V}$. The approximate solution $\rho^L$ is then to be found in $\mathcal{V}_L$, parametrized as $
    \rho^L(x) = \sum_{n=1}^{L}a_n\psi_n(x)$, and
to be tested against an arbitrary $\varphi\in\mathcal{V}_L$ of the form
$\varphi(x) = \sum_{m=1}^{L}k_m\psi_m(x)$.
In particular, we consider $k_m = e_m$ \color{black}to be the $m-$th standard unit vector\color{black}, leading to $\varphi = \psi_m$ for all $m=1,\dots,L$. This leads to the equation\begin{small}
\begin{equation*}\color{black}
\begin{aligned}
    -\beta^{-1}   \sum_{n=1}^{L} a_n \int\limits_{\Omega} \big(\nabla\psi_n+\psi_n\,\nabla V\big)\cdot\nabla\psi_m\,dx
    - \sum_{n,k=1}^{L} a_n a_k \int\limits_{\Omega} \big(\psi_n\big[(\nabla W*\psi_k)\big]\big)\cdot\nabla\psi_m\,dx = 0\,
\end{aligned}
\end{equation*}    
\end{small}
holding for all $m=1,\dots,L$. This system of nonlinear equation can be written in vector form for $\mathbf{a}:=\big(a_1,\dots,a_L\big)^\top$ as
\begin{equation}\label{nleq}
    0 = -(\beta^{-1} A + C)\mathbf{a} - b(\mathbf{a})
\end{equation}
with stiffness matrix $A\in\R^{L\times L}$, forcing term $C\in\R^{L\times L}$, and bilinear form $b(\cdot)\in\R^{L}$ defined as 
\begin{equation}\label{operators}
    \big[A\big]_{n,m} = \int\limits_{\Omega} \nabla\psi_n(x)\cdot\nabla\psi_m(x)\,dx\,, \qquad \big[C\big]_{n,m} = \int\limits_{\Omega} \psi_n(x)(\nabla V\cdot\nabla\psi_m(x))\,dx
\end{equation}
\begin{equation}\label{bil}
    b(\mathbf{a})=\begin{bmatrix}\color{black}
        \mathbf{a}^\top B_1\; \mathbf{a}\\
        \vdots\\\color{black}
        \mathbf{a}^\top B_L\; \mathbf{a}
    \end{bmatrix}
    \qquad\;\R^{L\times L}\ni\big[B_k\big]_{n,m}=\int\limits_{\Omega} \psi_n \big( \nabla W*\psi_k\big)\cdot\nabla\psi_m\,dx
\end{equation}

In the numerical tests, all the above integrals, including the convolution in \eqref{bil}, are computed by means of Gaussian quadrature rule. The number of quadrature points considered will be specified in each example.
Without loss of generality, we consider $\psi_1,\dots,\psi_L$ to be the first $L$ Fourier modes in $\Omega$, but any orthonormal basis of functions globally defined in space can be used instead. The main burden it terms of computational time (and memory) relies in the calculation and storage of the three-way tensor $B\in\R^{L\times L \times L}$, and can be partially mitigated by parallelization.

\subsubsection{Imposing the positivity and density conditions}
Since $\rho^L\approx\rho$ is the approximation of a probability density function we need to impose that it is strictly positive (for a strictly positive initial condition) and that it conserves mass  at all times, and in particular at the level of steady states 
\begin{equation*}
    \rho^L(x) >  0\quad\forall x\in\Omega\,,\qquad\qquad \int\limits_{\Omega} \rho^L(x)dx = 1\,.
\end{equation*}
The mass conservation constraint can be imposed directly on the system \eqref{nleq} by means of an additional matrix row. We sketch the idea in the two-dimensional case $x=(x^1,x^2)$, but it can be easily generalized to arbitrary dimensions. Assuming we have discretized the domain $\Omega$ in a mesh having $N+1$ nodes in both the dimensions, with the same discretization step $h$, we have
\begin{equation*}
   \sum_{i=1}^L a_i \sum_{j,k=1}^{N+1} \psi_i(x^1_{j},x^2_{k})h^2  \approx \int\limits_{\Omega} \rho(x)dx\,.
\end{equation*}
Moreover, using the notation $\zeta_i = \sum_{j,k=1}^{N+1} \psi_i(x^1_{j},x^2_{k})h^2$, the imposition of the density condition is equivalent to require
\begin{equation*} \Bar{\zeta}\mathbf{a} = 1\,,\qquad \text{where }\;\;
    \Bar{\zeta}= \begin{bmatrix}
        \zeta_1&\dots&\zeta_L
    \end{bmatrix}.
\end{equation*}
Then, our modified system of equations will be 
\begin{equation*}\label{neq_tilda}
    0 =  -\Tilde{A}\mathbf{a} - \Tilde{b}(\mathbf{a})\,.
\end{equation*}
with operators
\begin{equation*}
    \Tilde{A} = \begin{bmatrix}
        \beta^{-1}A\color{black}+C\\
        \Bar{\zeta}
    \end{bmatrix}
,\qquad
    \Tilde{b}(\mathbf{a}) = \begin{bmatrix}
        b(\mathbf{a})\\
        -1
    \end{bmatrix}.
\end{equation*}
On the other hand, the positivity contraint is handled by filtering out -- once a set of candidate steady states is generated -- all the $\rho_\infty$'s that display negative values for at least one point in the mesh.
%
%
\subsection{Deflation algorithm for PDEs}
After the spectral discretization and the imposition of the $L^1$ conservation constraint, we are left with the task of finding the (possibly) multiple solutions of the discretized Kirkwood-Monroe map 
\begin{equation}\label{residual_newton}\color{black}
    \mathcal{F}(\mathbf{a}) :=  -\Tilde{A}\mathbf{a} -\Tilde{b}(\mathbf{a}) = 0\,.
\end{equation}
In this section, we present the proposed numerical technique, based upon combining Newton's method with a deflation routine. We are motivated by root finding problems for scalar polynomials: if an iterative algorithm is used to find zeros $x_1,\dots,x_n$ of a polynomial $f(x)$, one can look for new, distinct roots of $f(x)$ by applying the same iterative algorithm to 
\begin{equation*}
    g(x) = \dfrac{f(x)}{\prod\limits_{i=1}^n (x-x_i)}\,.
\end{equation*}
The generalization of this method to the case of infinite-dimensional Banach spaces was proposed in \cite{deflation} and successfully applied in many different settings  \cite{deflation_2,deflation_3}. Here we briefly sketch the main elements of the deflation algorithm.

The procedure aims at systematically changing the PDE residual $\mathcal{F}$ so that known solutions are progressively ruled out. Namely, we are going to construct a new function $\mathcal{G}(\cdot)$ having all the zeros coinciding with the roots of $\mathcal{F}(\cdot)$, except for those we have already found. 
As a clarifying example, let us assume we have already identified a solution $\color{black}\mathbf{r}\color{black}$ of $\mathcal{F}(\mathbf{a}) = 0$, then
\begin{equation*}\label{defl0}
    \mathcal{G}(\mathbf{a}) = \mathcal{F}(\mathbf{a})\mathcal{M}(\mathbf{a},\color{black}\mathbf{r}\color{black})
\end{equation*}
where \color{black} for each $\mathbf{r}\in\Omega,\;\mathbf{a}\in \Omega\setminus\{\mathbf{r}\}$ \color{black} $\mathcal{M}(\mathbf{a},\color{black}\mathbf{r}\color{black})$ is a deflation operator, i.e. $\mathcal{M}$ is invertible, and $\mathcal{F}(\mathbf{a})\mathcal{M}(\mathbf{a},\color{black}\mathbf{r}\color{black}) = 0$ if and only if $\mathcal{F}(\mathbf{a})=0$.

Different classes of alternative choices of deflation operator $\mathcal{M}$ are proposed in \cite{deflation}, among which we consider the \emph{shifted exponential-norm deflation operator} that specifies 
\begin{equation}\label{defl}
    \mathcal{G}(\mathbf{a}) = \mathcal{F}(\mathbf{a})\bigg(\dfrac{\mathcal{I}}{\eta(\mathbf{a})} + \mathcal{I}\xi\bigg)\,,\qquad \eta(\mathbf{a}) = \|\mathbf{a}-\color{black}\mathbf{r}\color{black}\|^p\,,
\end{equation}
where $\mathcal{I}$ denotes the identity operator, $p$ is the deflation power and $\xi\geq0$  is the shift parameter. The deflating term $\eta(\mathbf{a})^{-1}$ ensures that $\mathcal{G}(\color{black}\mathbf{r}\color{black})\neq0$, while the shift operator $\mathcal{I}\xi$ forces the norm of the deflated residual to not vanish artificially as $\color{black}\|\mathbf{a}-\color{black}\mathbf{r}\color{black}\|\to\infty$. 

{\color{black}As discussed in \cite{deflation}, the selection of deflation parameters $p,\,\xi$ is specific to every problem. Here, as a general guideline, we consider $p > 1$ and $\xi > 0$. }
\subsubsection{Computing the derivative of the deflated operator}
In the following we assume that roots are found using Newton's method. Provided that, starting from an initial condition $\mathbf{a}_0$, the first root $r_{0}$ of $\mathcal{F}$ has already been found, the problem reduces to $\mathcal{G}(\cdot)$=0, for which the Newton's iteration reads
\begin{equation*}
\mathbf{a}^{n+1} = \mathbf{a}^n - \dfrac{\mathcal{G}(\mathbf{a}^n)}{\mathcal{G}'(\mathbf{a}^n)}\,.
\end{equation*}
This motivates a brief discussion of how to conveniently compute the derivative $\mathcal{G}'(\cdot)$. We focus on the case 
\begin{equation*}
    \mathcal{G}(\mathbf{a}) = \mathcal{F}(\mathbf{a})\bigg(\dfrac{\mathcal{I}}{\eta(\mathbf{a})}\bigg)\,,\qquad \eta(\mathbf{a}) = \prod_{i = 0}^{n} \|\mathbf{a}-r_{i}\|^p\,,
\end{equation*}
as the generalization to \eqref{defl} is straightforward. We write the deflation routine in a recursive way w.r.t. the number of roots. 

\begin{itemize}
    \item At the very first step, with only one known solution $r_{0}$ to $\mathcal{F}(\mathbf{a})=0$, we have
    \begin{equation*}
        \mathcal{G}(\mathbf{a}) = \mathcal{F}(\mathbf{a})\bigg(\dfrac{\mathcal{I}}{\eta_0(\mathbf{a})}\bigg)\,,\qquad \eta_0(\mathbf{a}) := \|\mathbf{a}-r_{0}\|^p\,,
    \end{equation*}
    with derivative
    \begin{equation*}
        \mathcal{G}'(\mathbf{a}) = \dfrac{\mathcal{F}'(\mathbf{a})}{\eta_{0}(\mathbf{a})} - \dfrac{\mathcal{F}(\mathbf{a})}{(\eta_{0}(\mathbf{a}))^2}\otimes\eta'_{0}(\mathbf{a})\,.
    \end{equation*}
    where $\otimes$ denotes the outer product.
    Assume that, by applying Newton's method for $\mathcal{G}(\mathbf{a})=0$ from the same initial guess $\mathbf{a}_0$, we find the second root $r_1$. 
    \item At the $n-$th iteration of the deflation loop, assuming to know roots $r_0,\dots,r_{n-1}$, the adjusted residual operator $\mathcal{G}$ reads
    \begin{equation*}
        \mathcal{G}(\mathbf{a}) = \mathcal{F}(\mathbf{a})\bigg(\dfrac{\mathcal{I}}{\eta_n(\mathbf{a})}\bigg)\,,\qquad \eta_n(\mathbf{a}) = \eta_{n-1}(\mathbf{a})\,\tau(\mathbf{a})\,,
    \end{equation*}
    for the temporal variable $\tau(\mathbf{a}) =\|\mathbf{a}-r_{n-1}\|^p$ and, dropping all the dependencies w.r.t. $\mathbf{a}$, we have the derivative
    \begin{equation*}
        \mathcal{G}' = \dfrac{\mathcal{F}'}{\eta_{n}} - \dfrac{\mathcal{F}}{\eta_{n}^2}\otimes\big(\eta_{n-1}'\,\tau + \eta_{n-1}\,\tau'\big)\,.
    \end{equation*}
    Newton's method is then used to find the next root $r_n$ of $\mathcal{G}(\mathbf{a})=0$.
\end{itemize}
The deflation loop stops when the root finding technique diverges for a given deflated operator and initial condition $\hat{\mathbf{a}}$. In the numerical tests we adopt the following stopping criteria for the Newton's iteration:  at the $n-$th step, either we have $\|\mathbf{a}^n-\mathbf{a}^{n-1}\|\leq\varepsilon$, with $\varepsilon=10^{-10}$, or $n\geq10^3$.

%
\subsection{The Haken-Kelso-Bunz Model}
In this section, we apply the root-finding procedure described above to a one-dimensional mean-field model on $[0, 2 \pi]$ with periodic boundary conditions. In particular, we consider the Haken-Kelso-Bunz (HKB) model~\cite[Sec. 5.3.3]{FrankFP} which was introduced as a model for human hand movement. This model exhibits phase transitions and has a rich set of (stable and unstable) stationary states; see also~\cite{angeli2023well}. We then compare the results of our numerical simulations with known analytical results.

\begin{figure}[!h]
    \centering
    \includegraphics[width=0.34\textwidth]{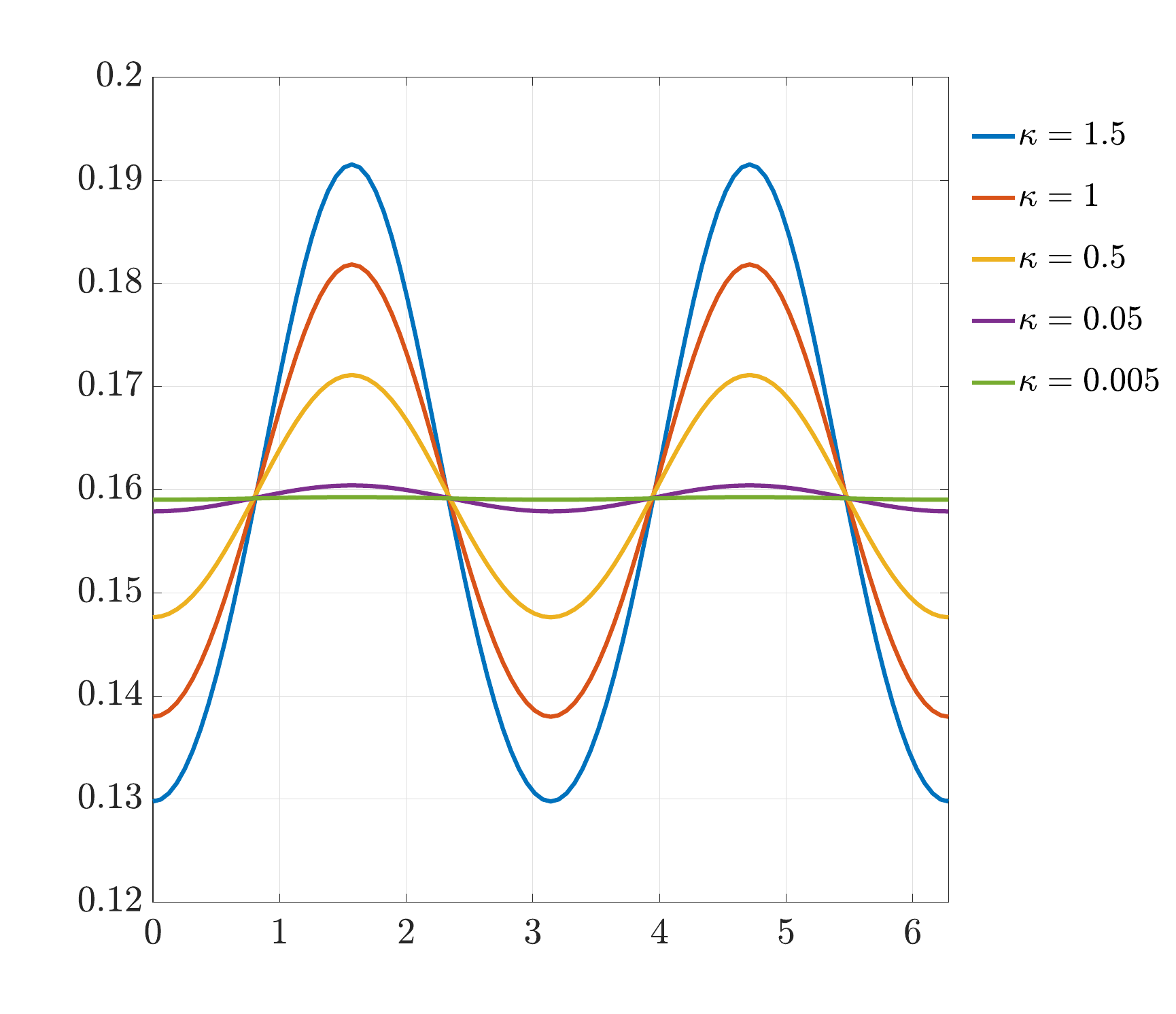}
    \includegraphics[width=0.295\textwidth]{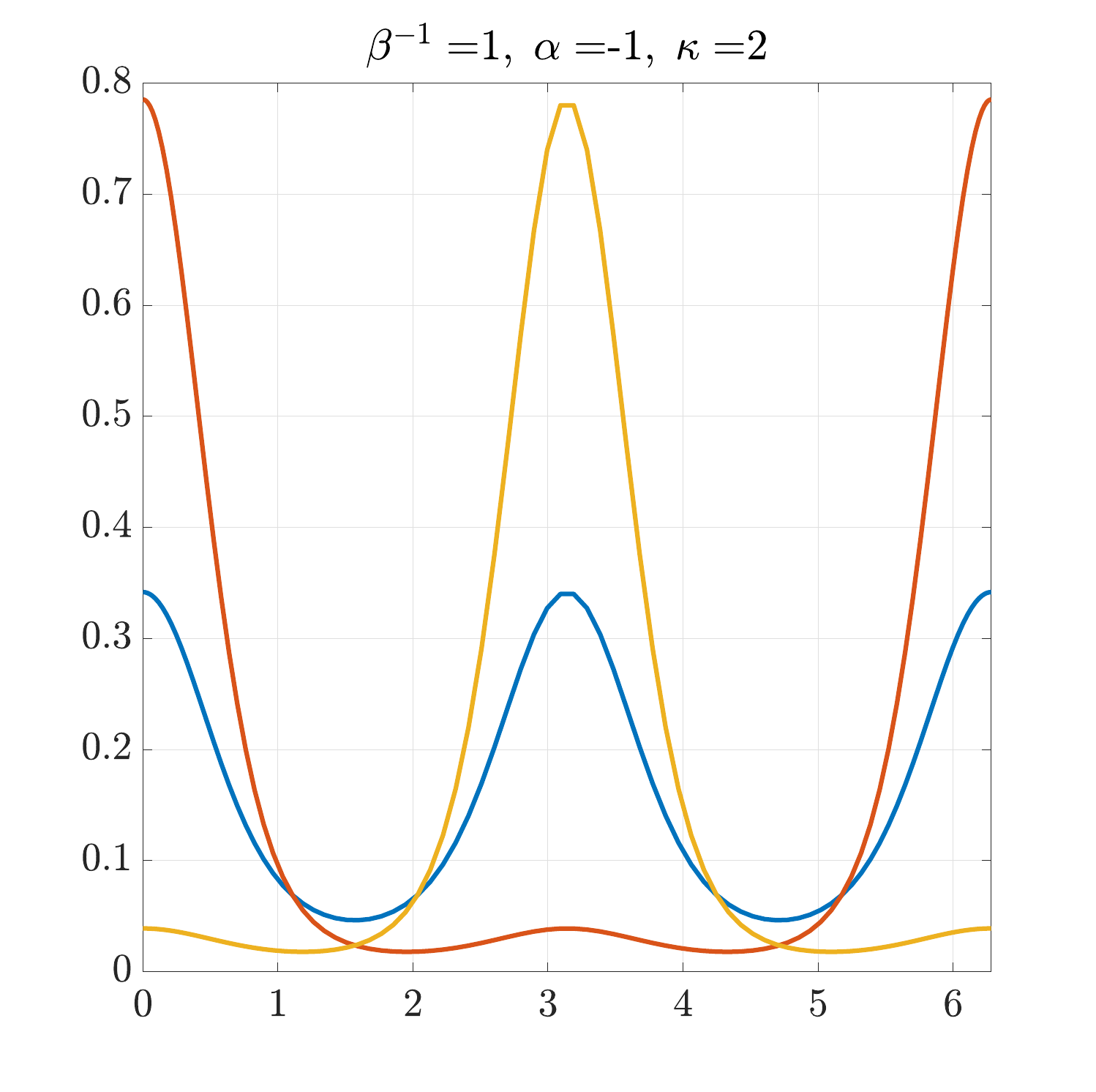}\hfill
    \includegraphics[width=0.295\textwidth]{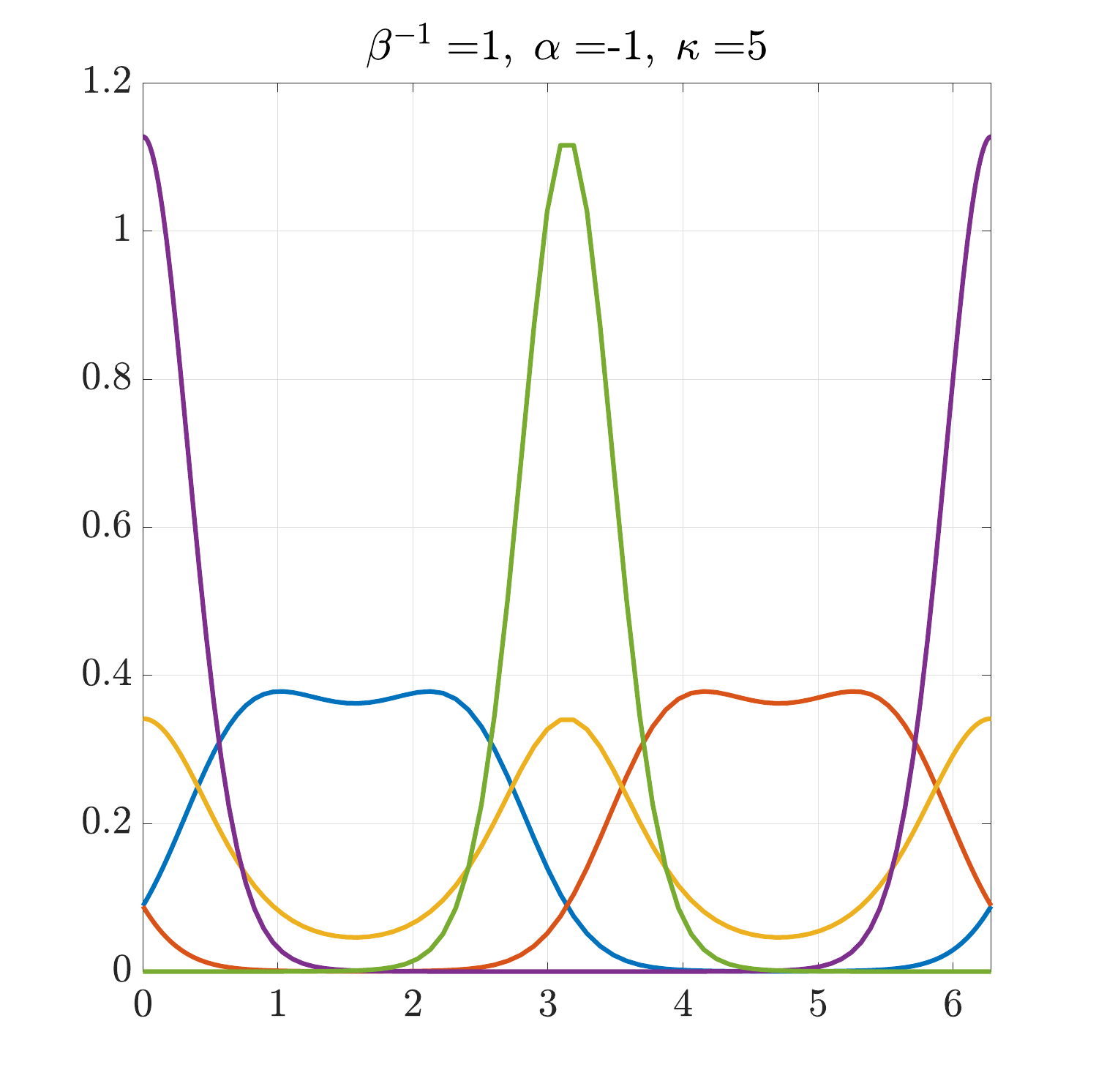}
    \caption{Steady states for the McKean–Vlasov PDE with potentials in Eq. \eqref{HKB_potentials} via deflated Newton's method for different interaction strength $\kappa$. Weak interaction regimes lead to a unique stable steady state, that collapses into the uniform distribution as $\kappa\rightarrow0$ (left). However, for stronger interactions with $\kappa = 2$ and $\kappa=5$ (centre and right, respectively), multiple stable and unstable stationary configurations coexist.}
    \label{fig:HKB1}
\end{figure}

We consider the PDE~\eqref{fp1} for $\rho:[0,2\pi]\to\mathbb{R}$ for confining and interaction potentials, respectively,
\begin{equation}\label{HKB_potentials}
\begin{aligned}
    W(x&,y)= -\kappa \,cos(x-y)\\
    V_{HKB}(x)& = -c_1\, cos(x) - c_2\, cos(2x)
\end{aligned}
\end{equation} 
with parameters $c_1, \, c_2 \in \mathbb{R}$ and $\kappa\geq0$. We consider the case for $c_1=0, c_2 = -\alpha$, $\alpha=-1$. For the initial guess in Newton's method, we consider a uniform random vector $a_0$ with entries normalized such that $\rho(0)$ is a distribution. We consider bifurcations with respect to the interaction parameter $\kappa$, and we fix the inverse temperature $\beta^{-1} = 1$. {\color{black} The choice of the model parameter values is motivated by the analysis provided in ~\cite[Sec. 5.3.3]{FrankFP}. } As expected, we see that, when the interaction strength $\kappa$ is below the critical value $\kappa<\kappa_c$, the PDE has a unique stationary solution. This stationary state, for the interaction and confining potentials~\eqref{HKB_potentials}, is calculated in~\cite{angeli2023well,bertoli2024}. 
However, when $\kappa$ is above the critical value $\kappa>\kappa_c$, we expect to have multiple steady states, as displayed in Figure \ref{fig:HKB1}.
\subsubsection{Verification of Solutions}
We solve the time-dependent McKean–Vlasov PDE~\eqref{fp1} for the HKB potential~\eqref{HKB_potentials}. The stability and bifurcation analysis presented in \cite[Sec. 5.3.3]{FrankFP}, on which we rely to validate the deflated Newton's method as a steady-state finding procedure.  

The McKean–Vlasov PDE with potential as in \eqref{HKB_potentials} reads
\begin{equation}\label{HKB_FP}
    \partial_t\rho = \beta^{-1} \partial_x^2 \rho +  \,\partial_x  \bigg[\rho\bigg(2\alpha\, sin(2x) + \kappa\int\limits_{\Omega} sin(x-y)\rho(y)dy\bigg)\bigg]\,,
\end{equation}
with steady states that are solutions of the Kirkwood-Monroe integral equation
\begin{equation*}
    \rho_\infty(x) = Z^{-1} \exp\bigg\{-\beta \bigg(\alpha\, cos(2x) - \kappa\int\limits_{\Omega} cos(x-y)\rho_\infty(y)dy\bigg)\bigg\},
\end{equation*}
for $Z$ normalization constant. Introducing the notation  
\begin{equation*}
    \big<g(X)\big> := \int\limits_{\Omega} g(x)\rho(t,x) dx\,,
\end{equation*}
we can write equation \eqref{HKB_FP} in the form
\begin{equation*}\color{black}
    \partial_t \rho = \beta^{-1} \partial_x^2 \rho +  \kappa\,\partial_x \bigg[\rho\bigg(2\alpha\, sin(2x)+\big(sin(x)\big<cos(X)\big> + cos(x)\big<sin(X)\big>\big)\bigg)\bigg]\,,
\end{equation*}
for which the steady states come as solutions of the integral equation 
\begin{equation}\label{rho_infty}
    \rho_\infty(x,\mathbf{m}) = Z(\mathbf{m})^{-1} \exp\bigg\{-\beta\bigg(\alpha\, cos(2x) -\kappa\big[cos(x)\,m_1 + sin(x)\,m_2\big]\bigg\}\,.
\end{equation}
Here, the parameters $\mathbf{m}=(m_1,m_2)$ are placeholders for the expectations $\big<cos(X)\big>$, $\big<sin(X)\big>$ w.r.t. $\rho_\infty$. This means we can solve equation \eqref{rho_infty} by coupling with the self-consistency equation 
\begin{equation}\label{self_consistency}
    \mathbf{m} = \color{black}\mathbf{r}\color{black}(\mathbf{m})\,\qquad \text{for }\; \color{black}\mathbf{r}\color{black} =\big(R_1,\;R_2\big)
\end{equation}

\begin{equation*}
    \text{for }\quad\color{black}\mathbf{r}\color{black}(\mathbf{m}) :=\bigg( \int\limits_{\Omega} cos(x)\rho_\infty(x,\mathbf{m})dx,\;\int\limits_{\Omega} sin(x)\rho_\infty(x,\mathbf{m})dx\bigg)\,.
\end{equation*}

\begin{remark}
Each one of the (possibly) many solutions $\mathbf{m}$ of the self-consistency equation \eqref{self_consistency} is uniquely linked to a steady state. Therefore, $\mathbf{m}$ serves as a vector-valued order parameter for detecting phase transitions.
Given a steady state $\Tilde{\rho}$, it is possible to estimate the associated 
 $\mathbf{m}$ via Linear Least Squares regression (LSS) for coefficients $a,b,c$ defined as
\begin{equation*}
     \underbrace{log(\Tilde{\rho})}_{Y} = - \beta\,V(x)  \underbrace{-log(Z)}_{a} +\underbrace{\beta\,\kappa\,cos(x)\,m_1}_{\phi_1(X)\,b} +\underbrace{\beta\,\kappa\,sin(x)\,m_2}_{\phi_2(X)\,c}\,.
\end{equation*}
We will use this to estimate the order parameters $\mathbf{m}$ associated with the stationary solutions that we find via deflated Newton's method, and then check if they are consistent with \eqref{self_consistency}. 
From the stability analysis developed in \cite{FrankFP}, we expect a phase transition between a monostable and a multistable regime, which is studied here separately with and without additional symmetry assumptions that are discussed below.
\end{remark}
\subsubsection{Symmetric solutions} 
\begin{figure}
    \centering
    \includegraphics[width = 0.3\textwidth]{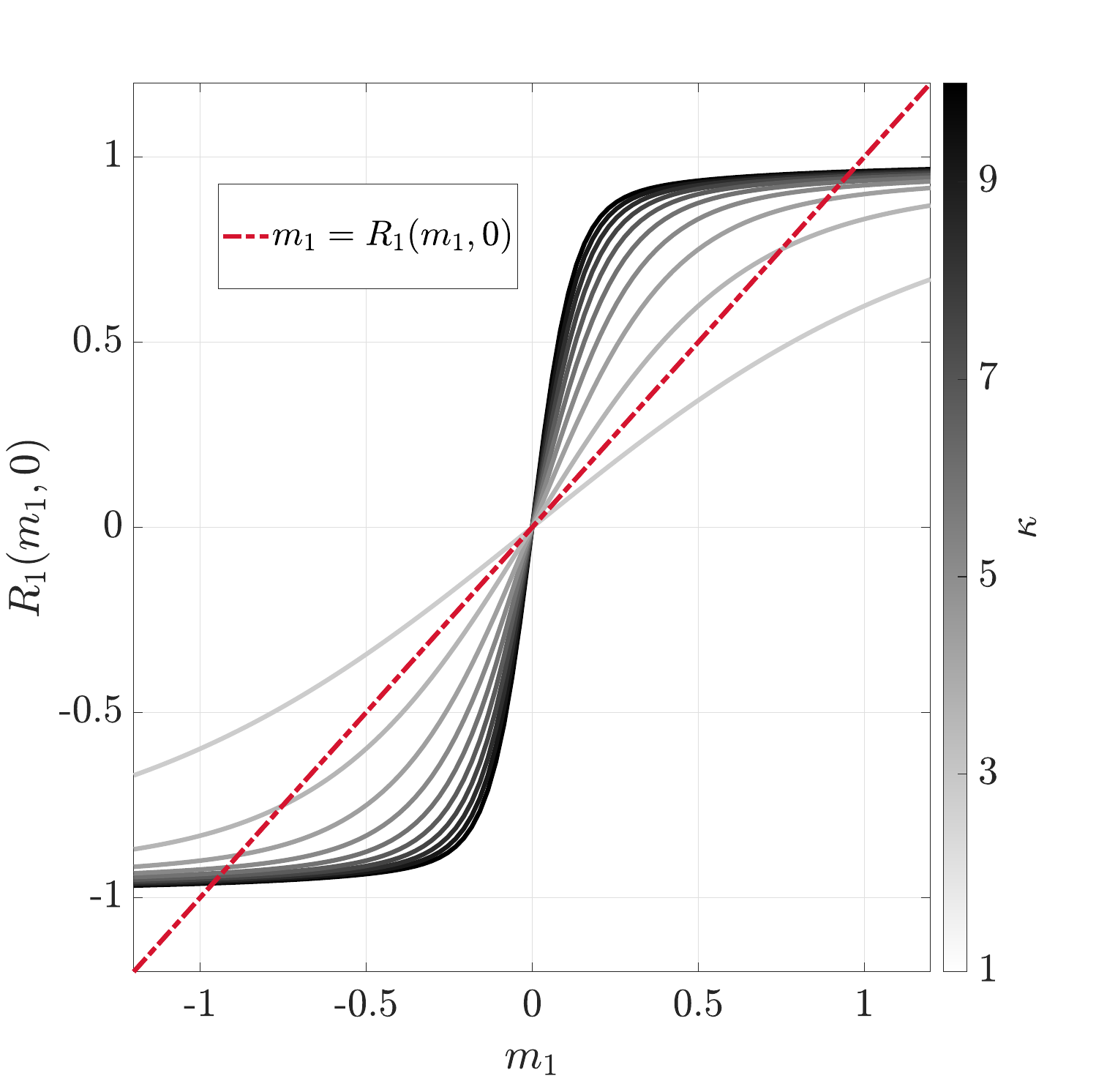}
    \caption{Solution to the self-consistency equation \eqref{self_consistency} for fixed $\beta=1$ and $\alpha=-1$, with varying interaction strength parameter $\kappa$. The dashed red line represents self-consistency, and its intersections with the curves identify symmetric stationary solutions of Eq. \eqref{rho_infty}. 
    In weak interaction regimes (similar to high noise $\beta^{-1}\gg1$), all solutions collapse into a single stable solution, indicating a saddle-node bifurcation.}
    \label{fig:self_consistency_plot_sym}
\end{figure}
\begin{figure}[h]
    \centering
    \includegraphics[width=0.34\textwidth]{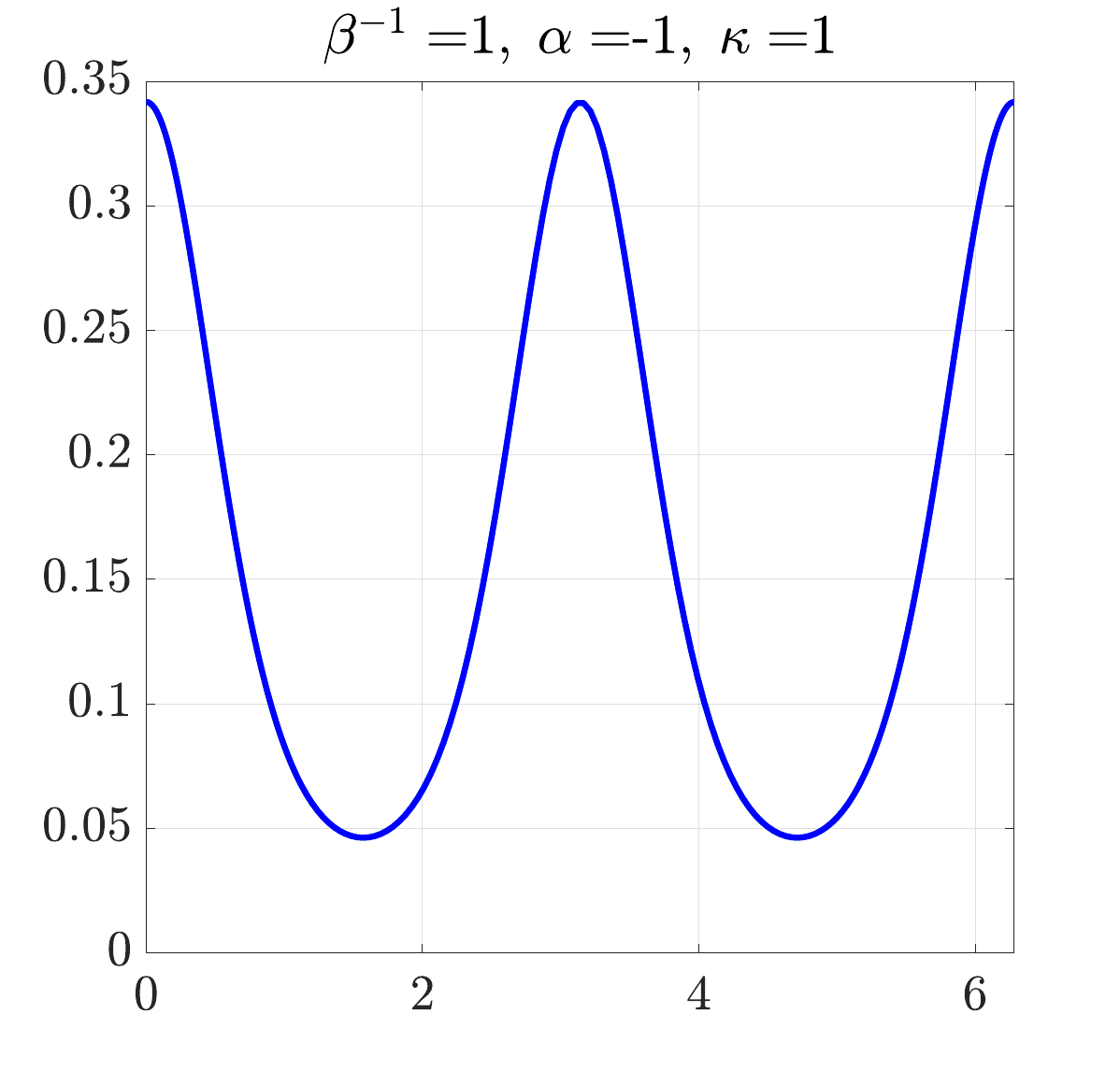}
    \includegraphics[width=0.34\textwidth]{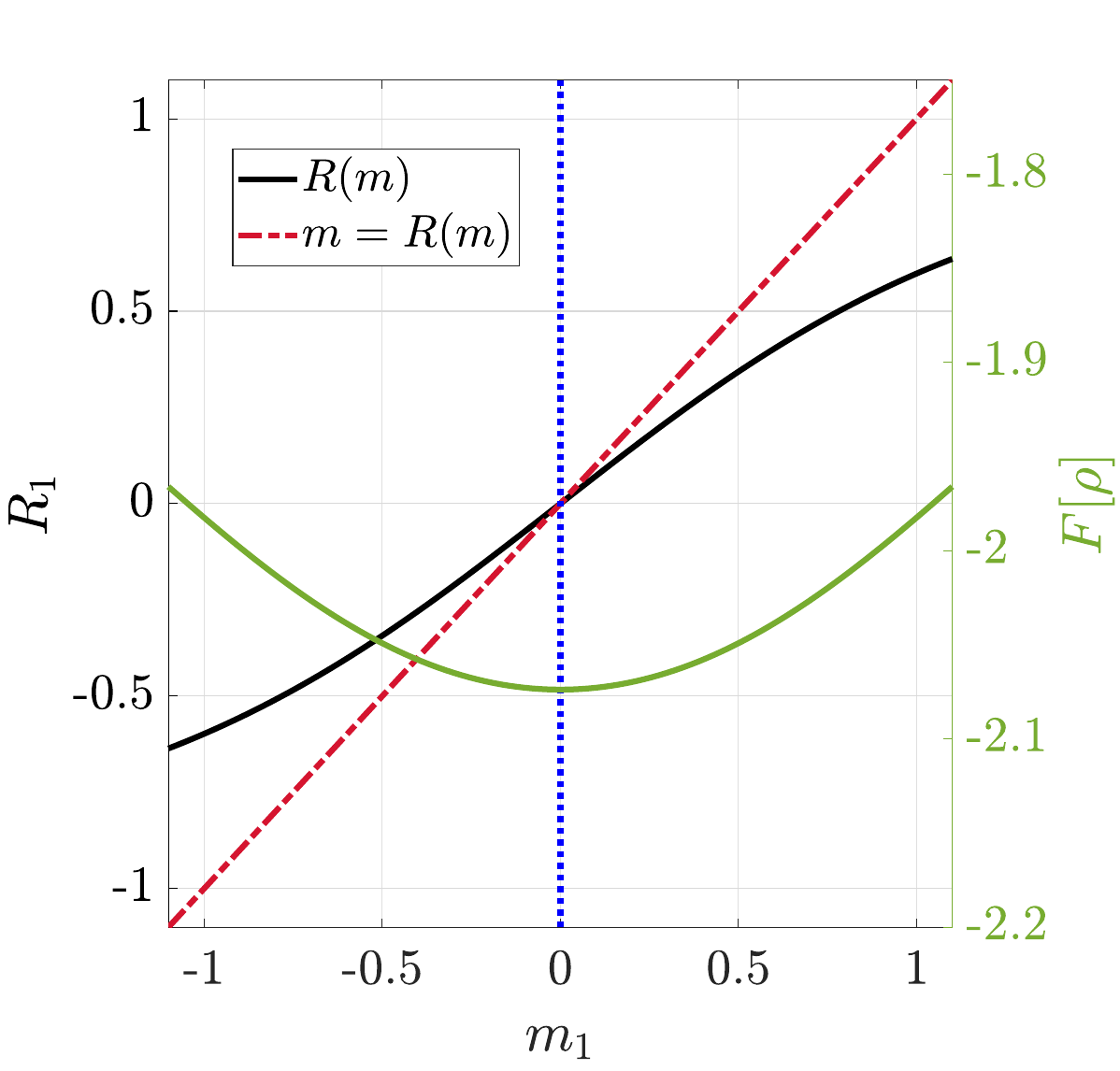}\\
    \includegraphics[width=0.34\textwidth]{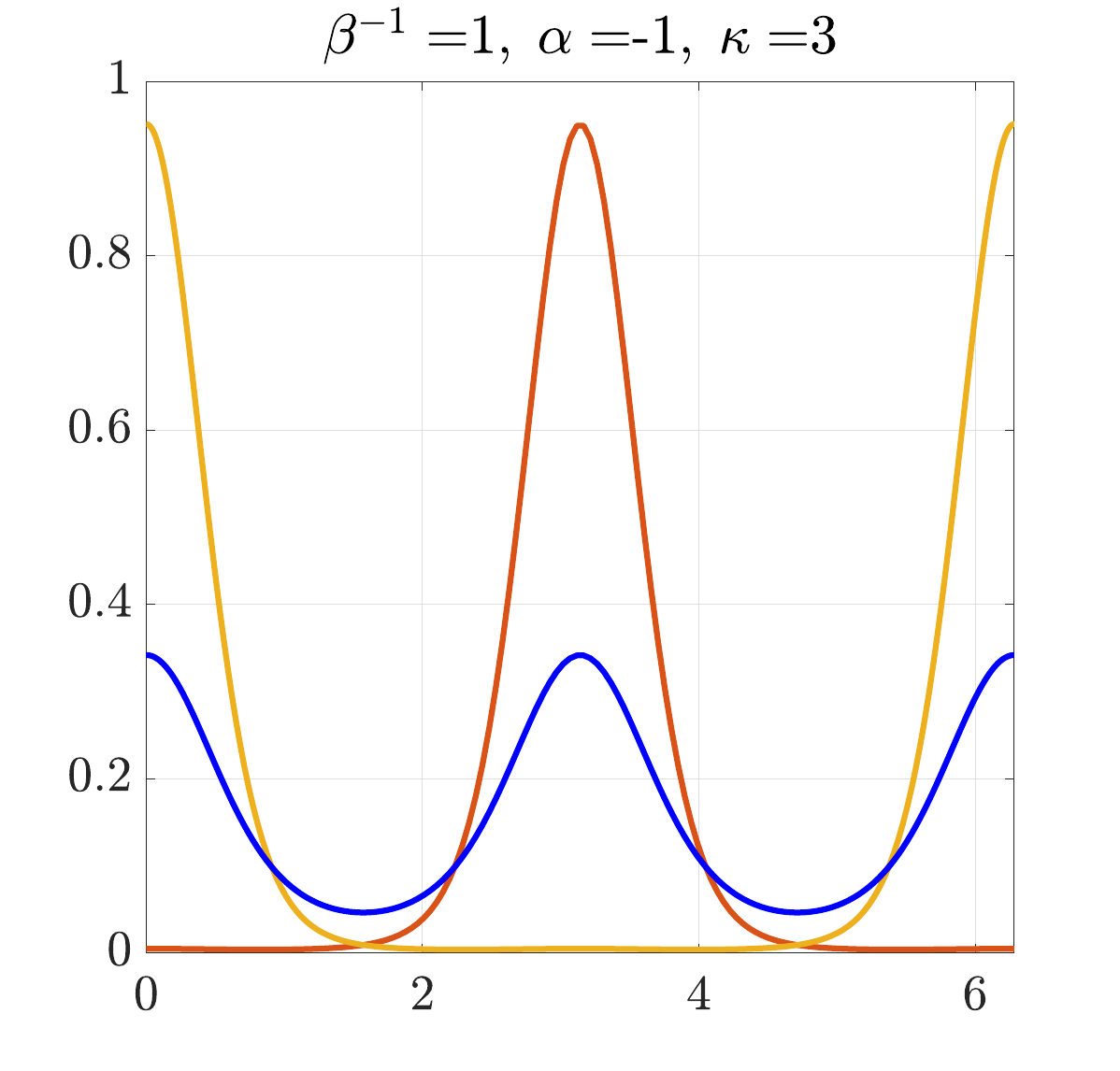}    
    \includegraphics[width=0.34\textwidth]{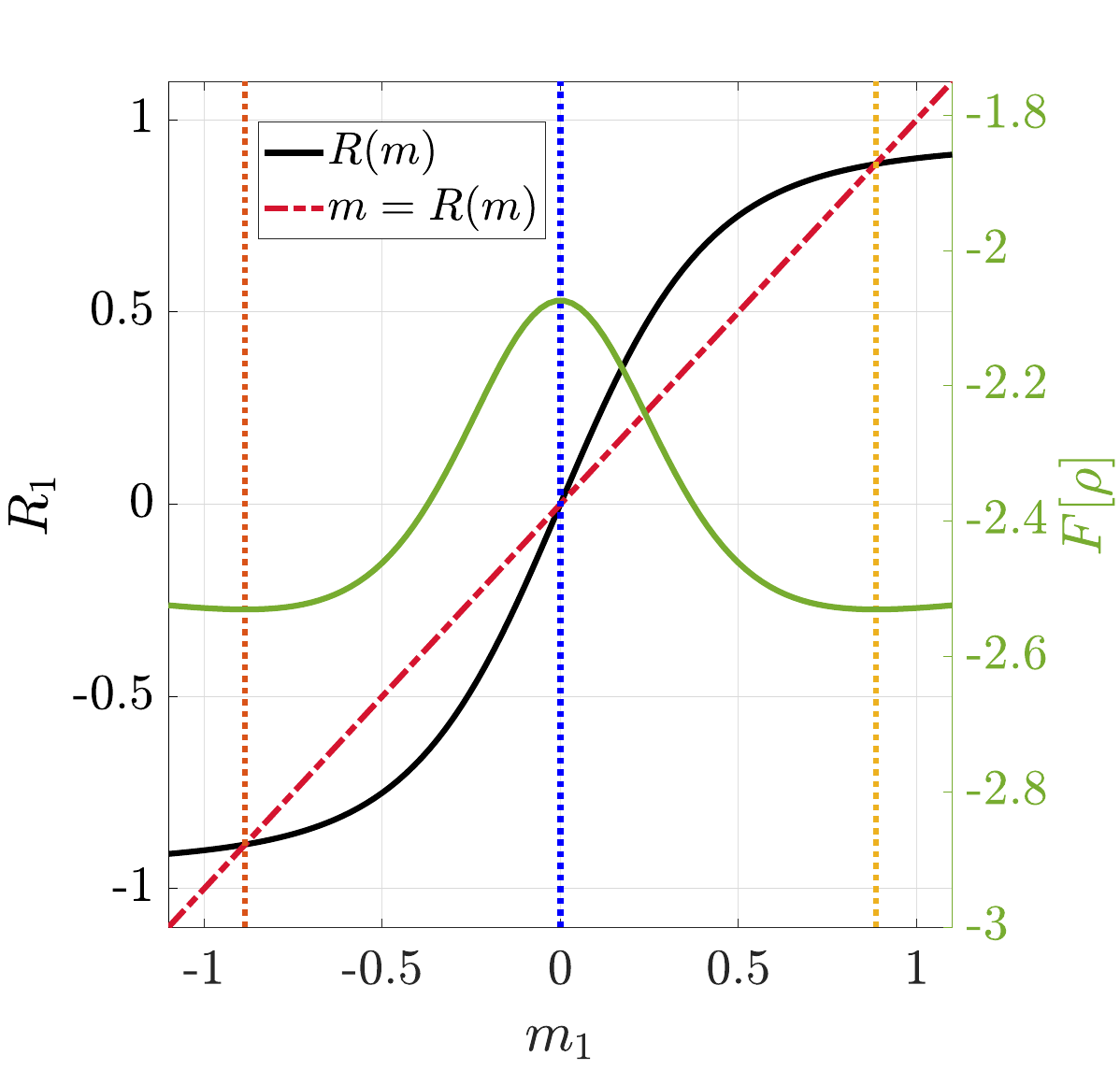}
    \caption{\textbf{Left}: Solutions obtained via the deflated Newton's method for two different parameter configurations, leading to a single solution (left) and triple solutions (right), respectively. \textbf{Right}: The self-consistency equation \eqref{sc_sym} and the free energy landscape. The dashed vertical lines represent the $m_1$ values estimated by the LSS method. The coloring of the curves corresponds to the solutions displayed on the left. All figures are generated with $\alpha=-1$, $\beta=1$, and varying interaction strength $\kappa=1$ (top) and $\kappa=3$ (bottom). The green curve shows the free energy $F[\rho]$ of the system, where stable and unstable configurations correspond to minimizing and non-minimizing critical points, respectively.
    }\label{fig:sc_symm}
\end{figure}

Due to the even symmetry of the confining potential $V(x)$, the second component $R_2(m_1, 0)$ vanishes. We focus then on solutions of the self-consistency equation (Eq. \eqref{self_consistency}) of the form:
\begin{equation}\label{sc_sym}
m_1 = R_1(m_1, 0), \qquad m_2 = 0,
\end{equation}
which correspond to symmetric steady state configurations. 
Figure \ref{fig:self_consistency_plot_sym} illustrates how the number of solutions to equation \eqref{sc_sym} varies with the inverse temperature parameter. Following \cite{FrankFP}, we expect the non-negative solutions $m_1 \geq 0$ to parameterize stable steady states, while negative solutions lead to unstable configurations. 
The self-consistency argument ensures that at least one symmetric stable steady state always exists (associated with $R_1(0,0) = R_2(0,0) = 0$). Figure \ref{fig:sc_symm} verifies that the deflation algorithm finds all the expected symmetric solutions. This is obtained by considering $l=100$ Fourier modes, $200$ Gauss-points, $p=2,\,\xi=1$, and with initial condition $\mathbf{a}^0$ given by the zero vector. All the steady states obtained via Newton's method give a residual \eqref{residual_newton} $\|\mathcal{F}(a)\|_{l^2}\leq 10^{-15}$.

\subsubsection{Dropping the symmetry assumption}

\begin{figure}[h]
        \centering
        \includegraphics[width=0.31\textwidth]{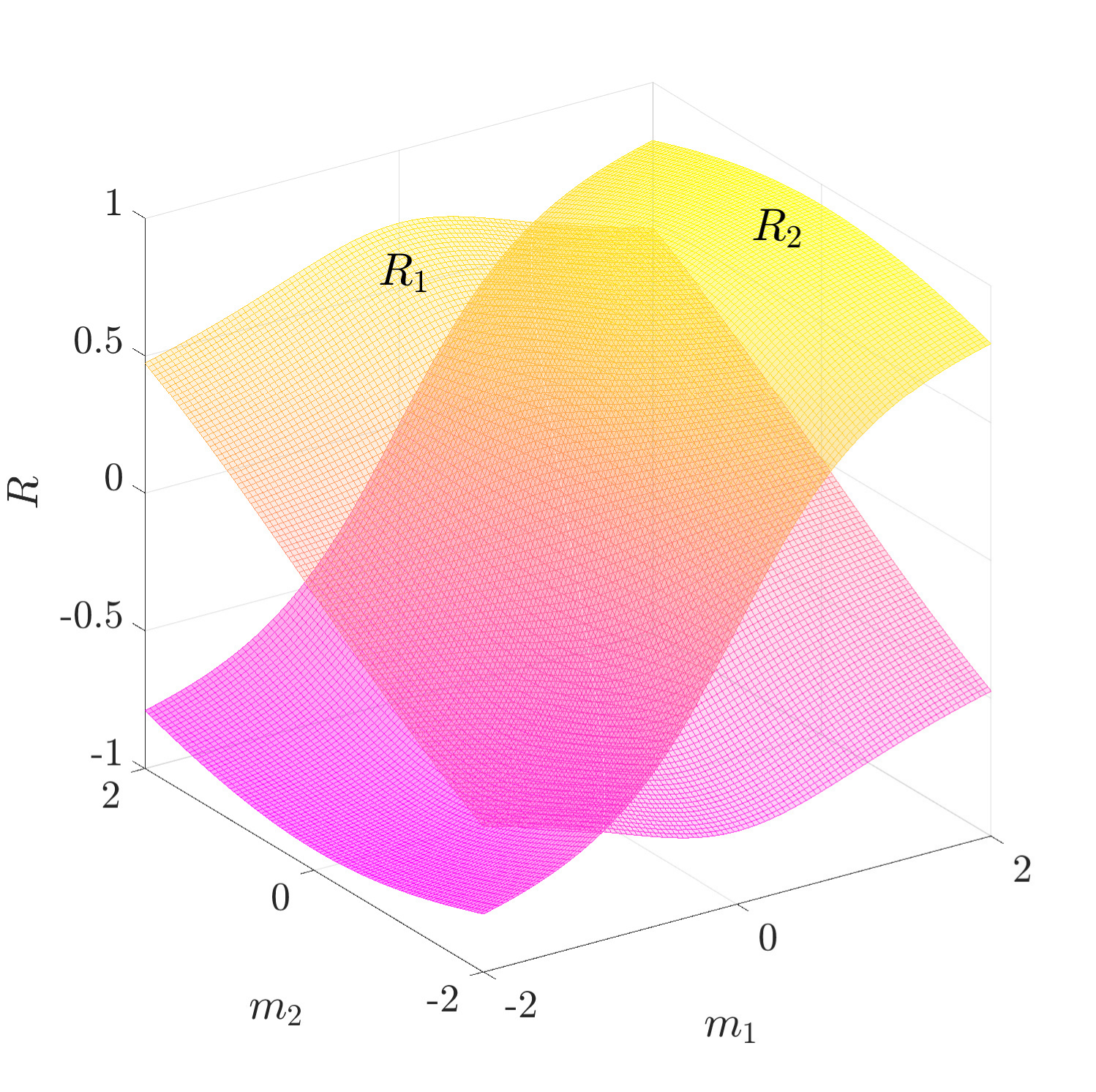}
        \includegraphics[width=0.65\textwidth]{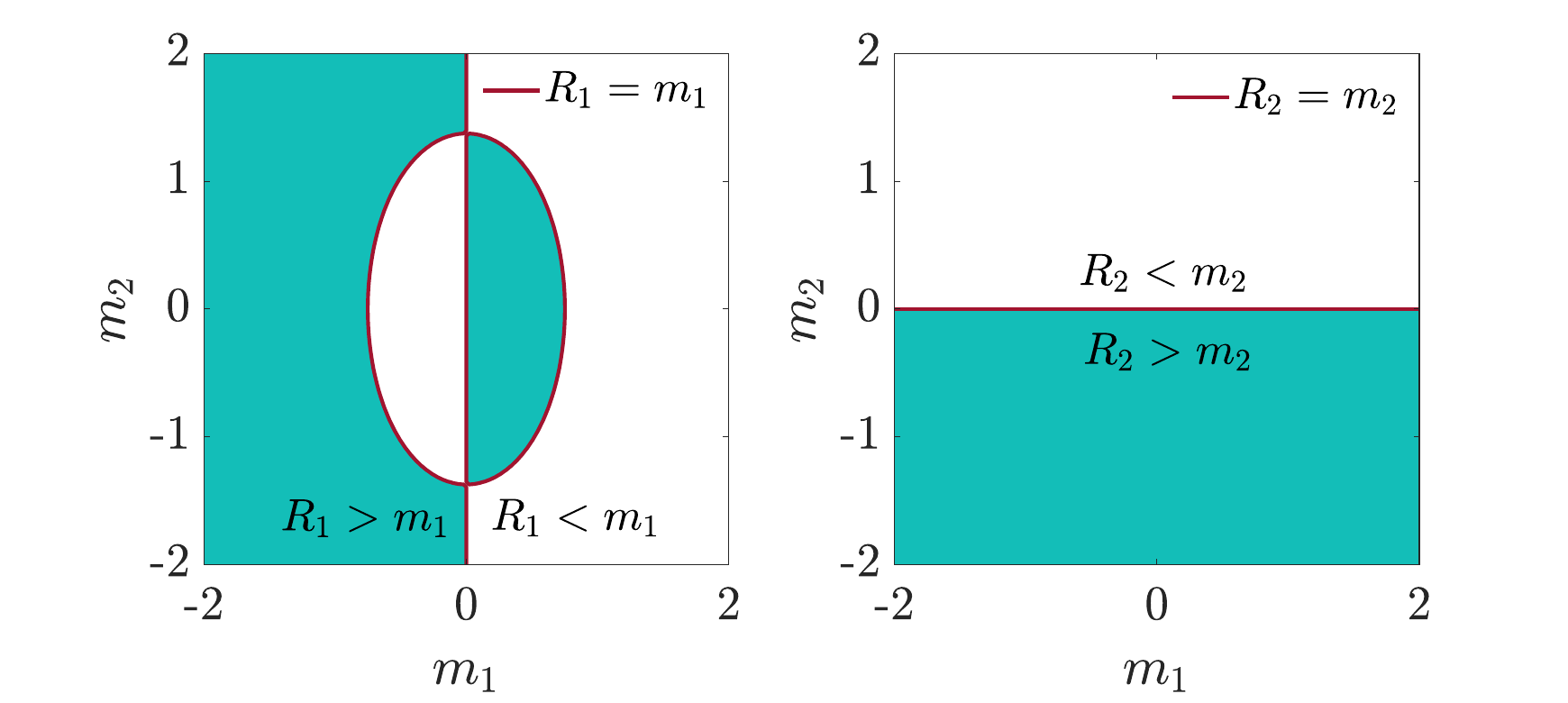}
        \caption{$\mathbf{\kappa=2}$: The complete self-consistency equation is linked to two surfaces, one for each component of $\color{black}\mathbf{r}\color{black}(\mathbf{m}
    )$ (left). The solutions $\mathbf{m}$ lie then on the diagonal planar section of those surfaces w.r.t. $m_1$ and $m_2$ respectively (right).}\label{fig:kappa_2}
\end{figure}

\begin{figure}[h]
    \centering
    \includegraphics[width=0.31\textwidth]{rho1m12} \hfill
    \includegraphics[width=0.31\textwidth]{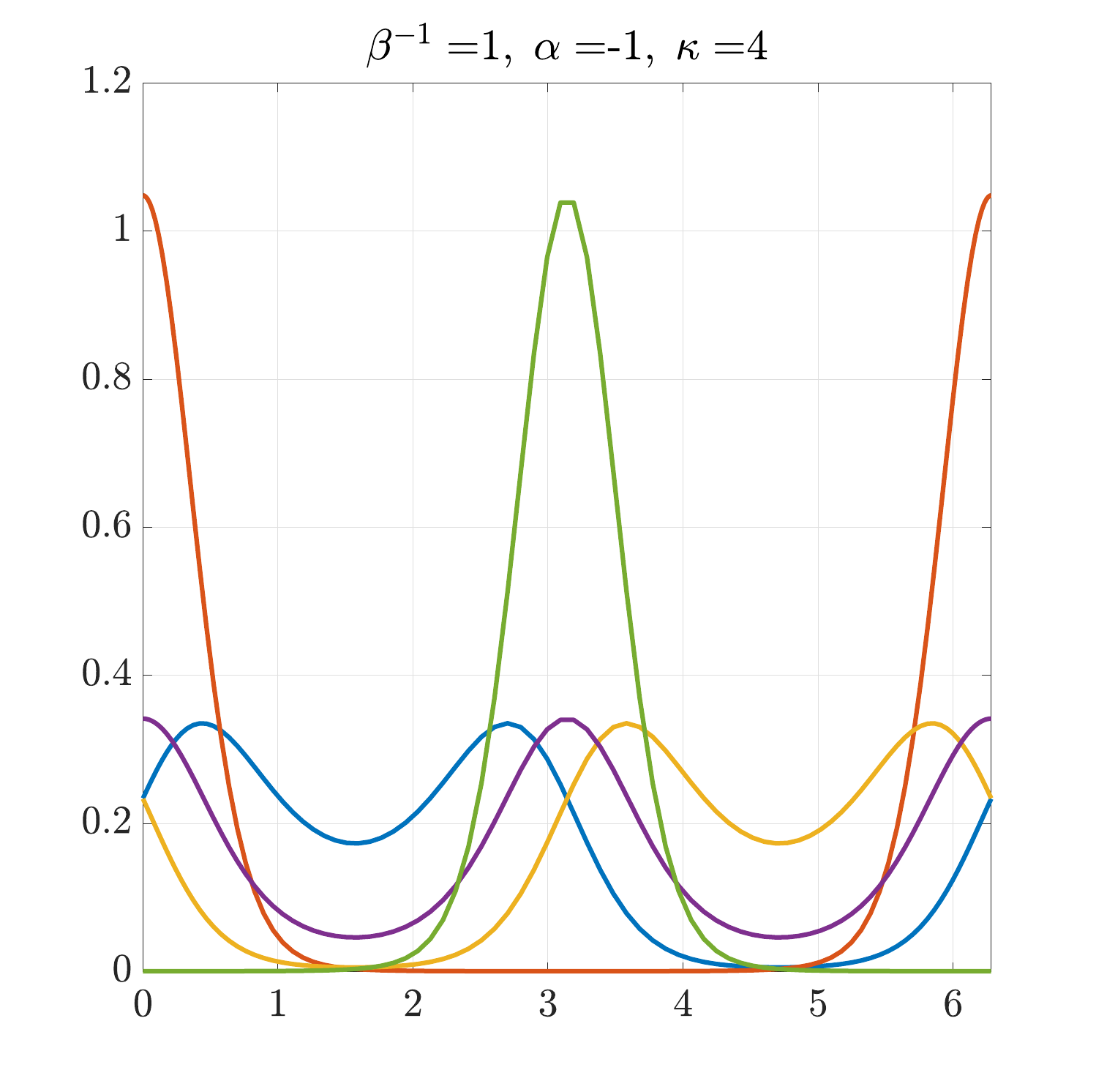}\hfill
    \includegraphics[width=0.31\textwidth]{rho1m15}\\
    \includegraphics[width=0.33\textwidth]{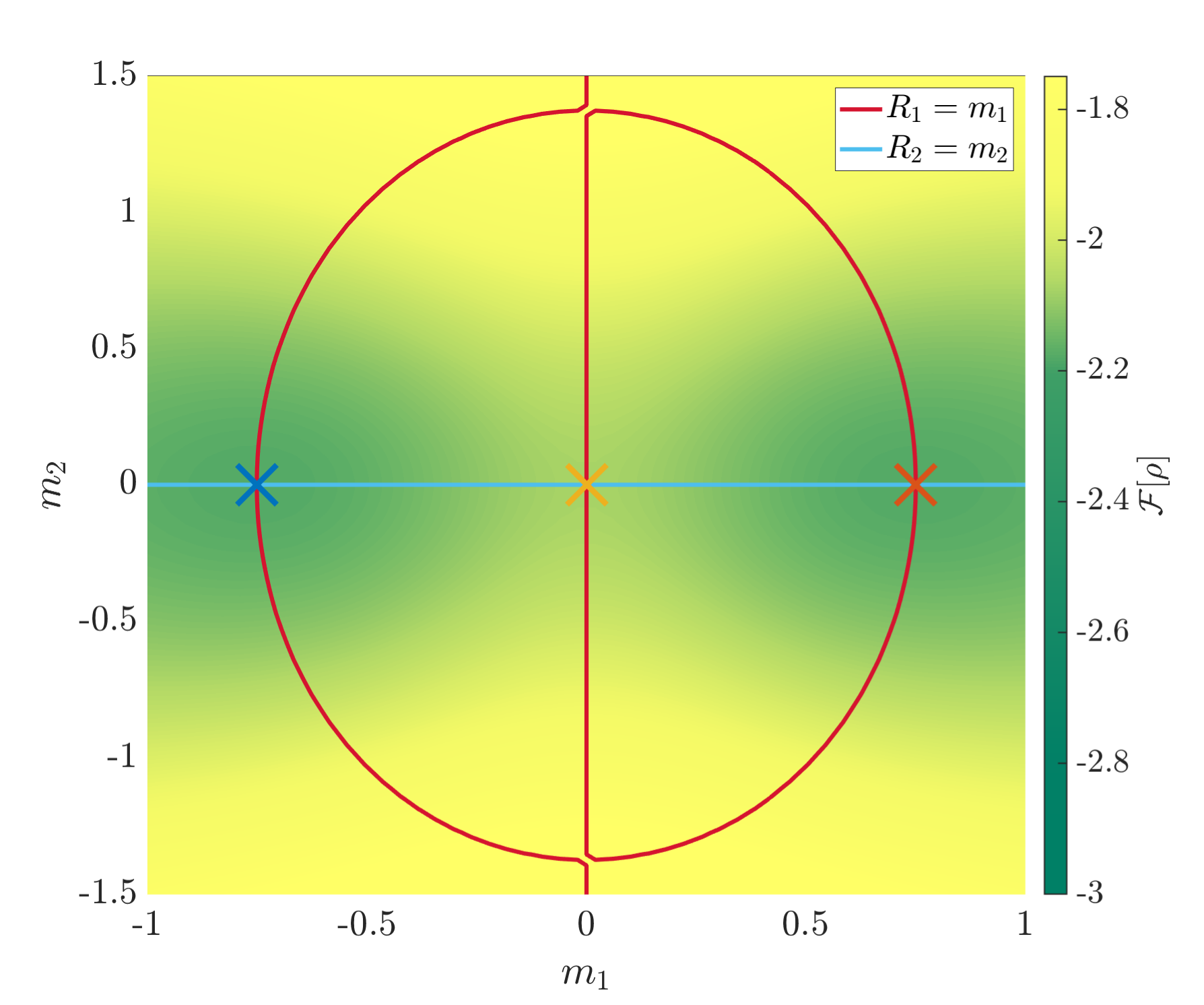} \hfill
    \includegraphics[width=0.33\textwidth]{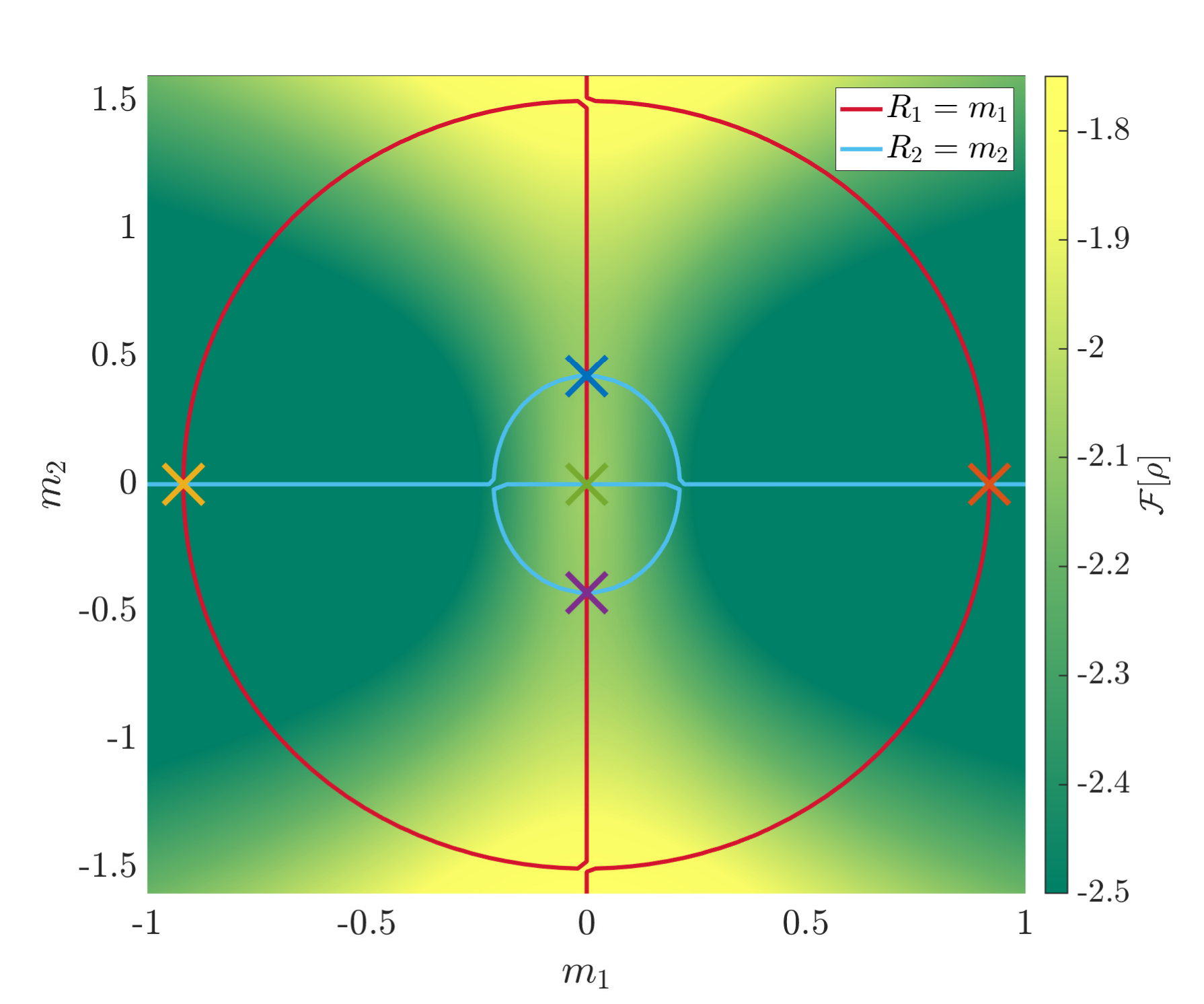}\hfill
    \includegraphics[width=0.33\textwidth]{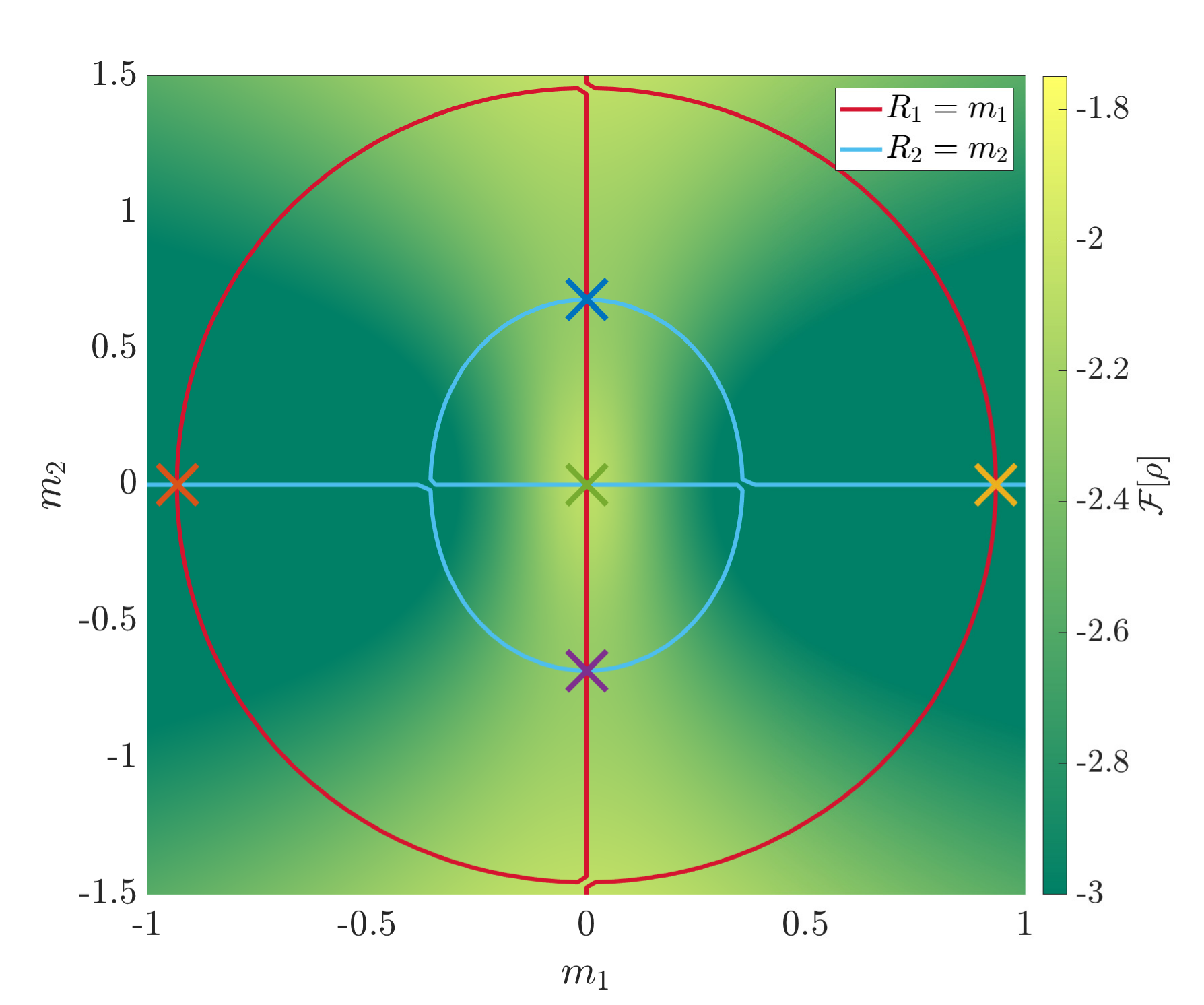}
    \caption{We display all the solutions (symmetric and asymmetric) found via deflation and we compared the estimated values of $m_1$ and $m_2$ with the contour plot and free energy landscape for different parameter configurations: $\mathbf{\kappa=2}$ (left), $\mathbf{\kappa=4}$ (center), $\mathbf{\kappa=5}$ (right). The estimated parameters lie in the desired intersections, validating the numerical solutions.
    }\label{fig:sc_asymm}
\end{figure}

For validating both symmetric and non-symmetric solutions (as in Figure \ref{fig:HKB1}), the full self-consistency equation \eqref{self_consistency} with  $\mathbf{m} = (m_1,\,m_2)$ is needed. We verify these solutions by comparing the LSS-estimated $\mathbf{m}$ with the intersection of the diagonal planar sections in Figure \ref{fig:kappa_2}. 
In Figure \ref{fig:sc_asymm} we superimpose the selfconsistency solutions over the free energy levels; this allows us to identify stable and unstable configurations. The figures have been obtained by considering $L=50$ Fourier modes, $100$ Gauss-points, $p=4,\,\xi=1$, and with initial condition $\mathbf{a}^0$ given by a random uniform vector normalized such that the associated $\rho$ is a distribution. All the steady states obtained via Newton's method give a residual \eqref{residual_newton} $\|\mathcal{F}(a)\|_{l^2}\leq 10^{-15}$.

%
%
\subsection{Phase transitions for the $O(2)$ model with magnetic field} 
\label{subsec:O2}

To validate the proposed root-finding procedure numerically, we apply it to the problem of the $O(2)$ model with an external magnetic field; for this problem, the stationary states can be computed analytically~ \cite[Lemma 1.26]{Degaldino}. We consider the state space $\Omega=[0,2\pi]$ with interaction and confining potentials
\begin{equation*}
    W(x,y) =  -cos(2\pi (x-y))\,,\qquad V(x)=-\eta cos(2\pi x)\,,
\end{equation*}
respectively, with $\eta\in(0,1)$. {\color{black} The presence of a small
confinement (e.g. $\eta=0.05$) is enough to break the translation invariance of the problem.} In~\cite[Lemma 1.26]{Degaldino}, it is proved that there exists a critical temperature $\beta^{-1}_c$ such that:
\begin{itemize}
    \item for $\beta^{-1}>\beta^{-1}_c$ there exists a unique stable steady state parameterized by 
\begin{equation}
    \rho_{min} = Z^{-1}_{min} e^{\alpha_{min}cos(2\pi x)}\,,\qquad Z_{min} = \int\limits_{\Omega} e^{\alpha_{min}cos(2\pi x)} dx\qquad\text{with $\alpha_{min}>0$;}
\end{equation}

\item for $\beta^{-1}<\beta^{-1}_c$ we expect \emph{at least} 2 steady state solutions of the form 
\begin{equation*}
\begin{aligned}
\rho_{min} &= Z^{-1}_{min} e^{\alpha_{min}cos(2\pi x)}\,,\qquad Z_{min} = \int\limits_{\Omega} e^{\alpha_{min}cos(2\pi x)} dx\\
        &\rho_{*} = Z^{-1}_{*} e^{\alpha_{*}cos(2\pi x)}\,,\qquad\;\;\; Z_{*} = \int\limits_{\Omega} e^{\alpha_{*}cos(2\pi x)} dx
\end{aligned}
\end{equation*}
for $\alpha_*<0<\alpha_{min}$ with both constants depending on $\beta^{-1}$. At low temperatures, $\rho_{min}$ is the unique minimizer ($\rho_*$ being a non-minimising critical point) of the free energy. This classifies $\rho_{min}$ and $\rho_*$ as stable and unstable steady states, respectively.
\end{itemize}
\begin{figure}[t]
    \centering
    \includegraphics[width=0.33\textwidth]{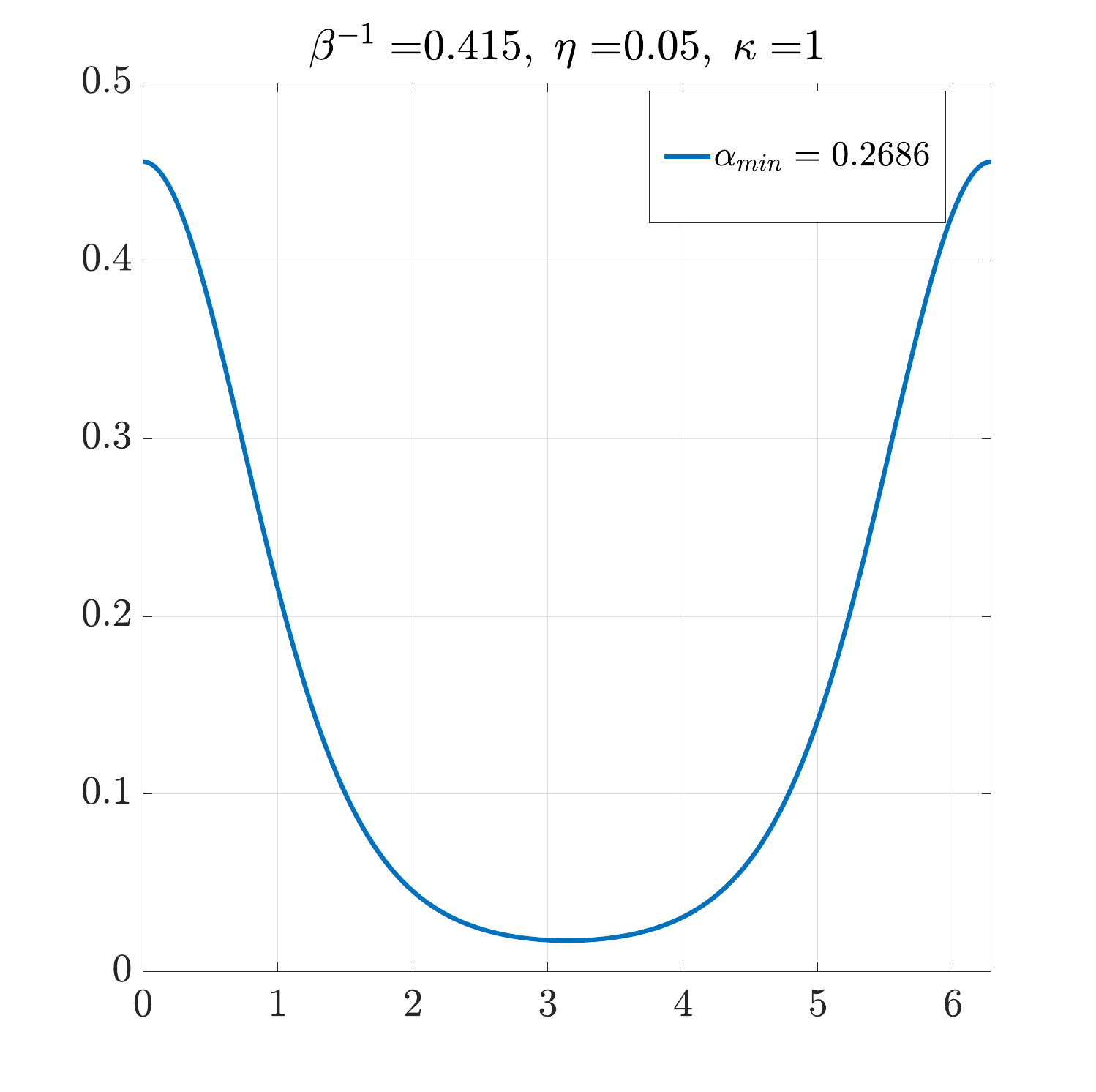}\hfill\includegraphics[width=0.33\textwidth]{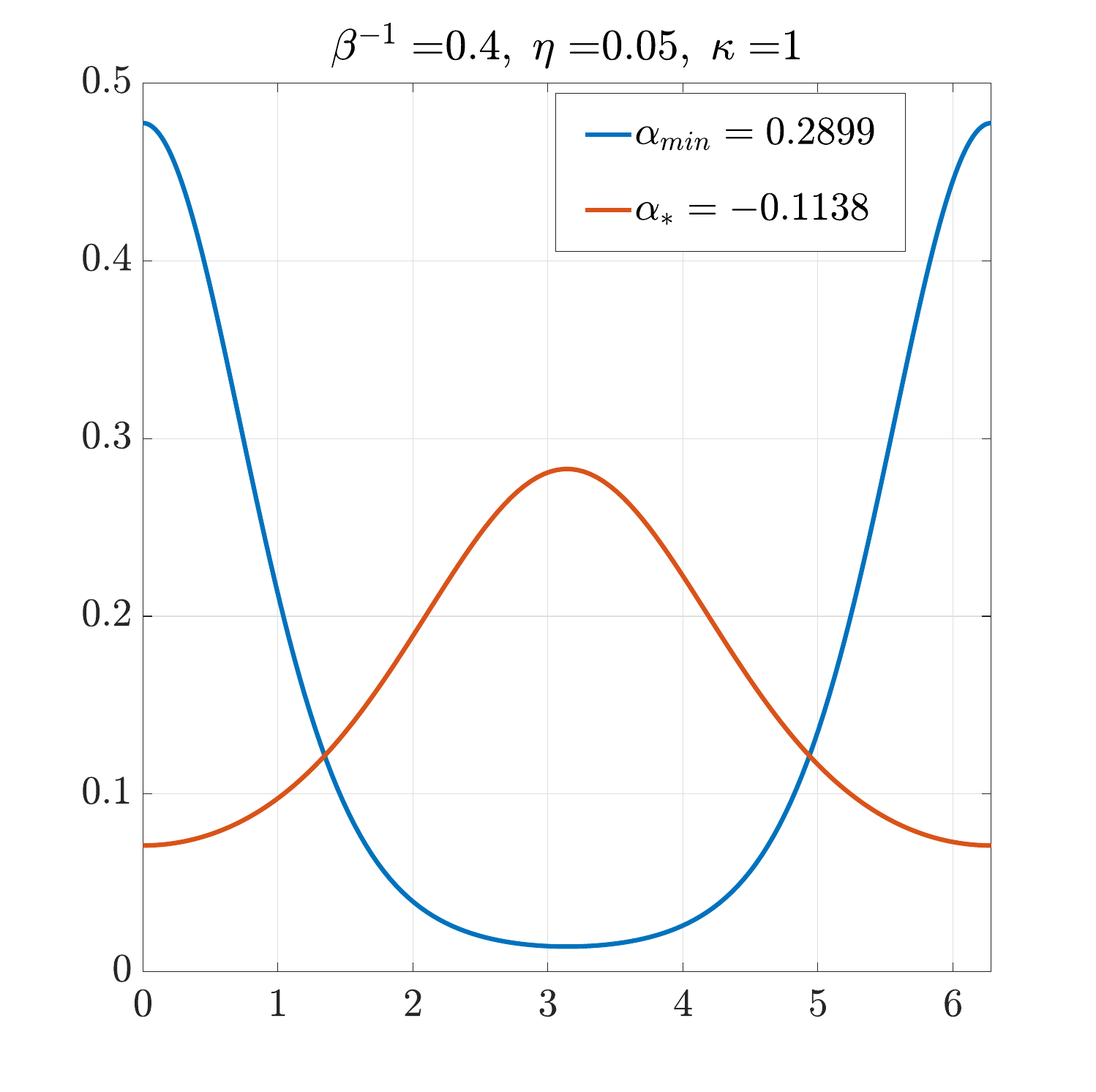}\hfill \includegraphics[width=0.33\textwidth]{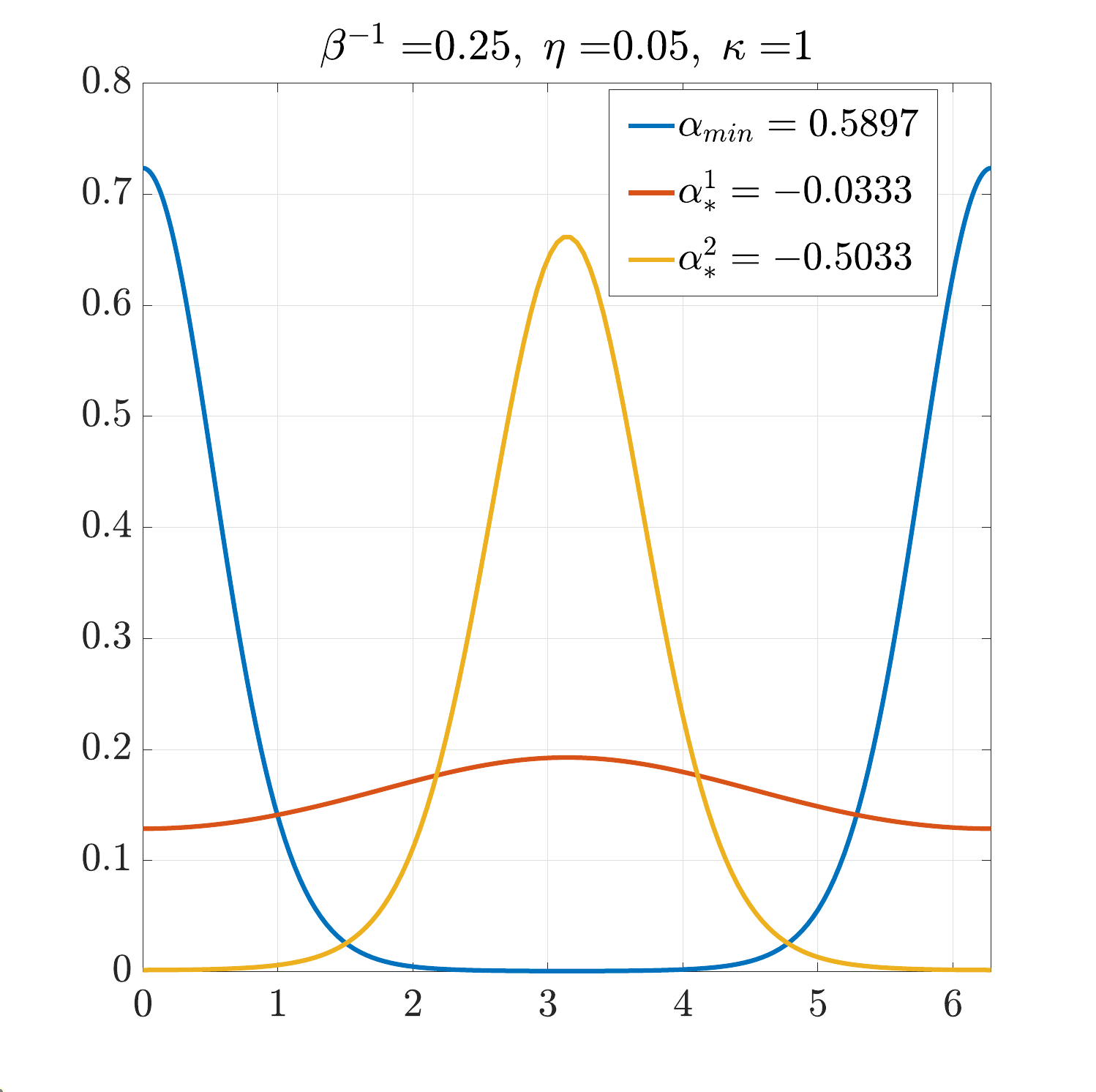} 
    \caption{Steady states for the modified noisy Kuramoto model (confining strength parameterized by $\eta = 0.05$) obtained via deflated Newton's method. We label the solutions with the LSS-estimated coefficients $\alpha_{min}$ and $\alpha_*$'s. }
    \label{fig:Lemma126}
\end{figure}
The results of our numerical experiments are displayed in Figure \ref{fig:Lemma126} \color{black}, and are obtained with the following parameters choice: $L=10$ Fourier modes, $100$ Gauss-points, $p=2$, $\xi=0.05$, and $ \mathbf{a}^0$ is the uniform distribution. \color{black} The coefficients $\alpha_{min}$ and $\alpha_*$'s, estimated using the LSS method as done in the previous example, match the expected behaviour. In particular, regardless the choice of $\beta^{-1}$, we always retrieve a positive $\alpha_{min}>0$, associated with a stable configuration. Below the critical temperature, which we place numerically between $0.4<\beta^{-1}_c<0.415$, other unstable steady states, parametrized by $\alpha_*<0$, start to appear.

%
%
\section{Optimal feedback stabilization of unstable steady states}
\label{sec:stabilization}
Having determined all possible steady state solutions to our problem, we now turn our attention to the design of control laws which can effectively steer the McKean-Vlasov PDE towards a given stationary distribution. For this, we consider the mean-field equation with a control term in the drift
\begin{equation}\label{ocp}
	\partial_t \rho = \beta^{-1}\Delta\rho + \nabla \cdot\bigg(\rho\big((\nabla W \ast \rho)+\nabla V + u\big)\bigg)\,,\qquad (t,x)\in\Omega\times[0,T]\,,
\end{equation}
where $\color{black} u=u(t,x)\in L_2(0,T;L_\infty(\Omega))$ is a time-dependent vector field representing an external control signal. \color{black} In the following we refer to $\rho\in \mathcal{W}(0,T)$ as the weak solution to \eqref{ocp}, where the Bochner space reads 
$$\mathcal{W}(0,T):= \{\rho\in L_2(0,t;H^1(\Omega))\,|\,\partial_t\rho\in L_2(0,T;H^{-1})\}\,,$$
with $H^{-1}$ dual space of $H^1$. \color{black}

A natural control problem arising in this context is the acceleration of convergence towards a stable steady state, which can be achieved by studying the linearized control problem around an equilibrium \cite{control2018breiten}. In this paper, our goal is to choose a control law  to steer the solution to the McKean-Vlasov PDE towards a prescribed unstable steady state. The synthesis of such a control signal is achieved by casting the following dynamic optimization problem:
\begin{equation}\label{ocp2}
\underset{u}{\min}\;\cJ(u):=\frac12\int_0^T\|\rho(t,\cdot)-\rho_{\infty}\|^2_{L_2(\Omega)}+\gamma\|u(t,\cdot)\|^2_{L_2(\Omega,\R^d)}\,dt\,+\Psi(\rho(t,x))\,,
\end{equation}
subject to \eqref{ocp}. In this formulation $T>0$ denotes a prescribed prediction horizon for stabilization, $\rho_{\infty}$ corresponds to the desired unstable steady state to be tracked, $\gamma>0$ is a control energy penalty, and $\Psi(\rho)$ is a terminal cost. In finite horizon optimal control problems, this terminal penalty plays a crucial role in enforcing stabilization towards the terminal state, representing a relaxation of the terminal constraint $\rho(t,x)=\rho_{\infty}$. In our case, we consider a terminal penalty given by
\begin{equation*}
    \Psi(\rho) = \eta\,\|\rho(t,\cdot)-\rho_\infty\|_{L_2(\Omega)}^2\,,\qquad \eta\geq0\,.
\end{equation*}
\begin{remark}\color{black} The non-negativity and mass conservation conditions have not yet been addressed in the optimal control problem. However, these properties are inherently satisfied when the initial condition $\rho_0\in L_2(\Omega)$ is non-negative and has an integral equal to $1$. For a more detailed discussion on this topic, we direct the reader to Corollary 2.4 in \cite{Borzi_secondorder}.
\end{remark}
%
%
\subsection{First-order necessary optimality conditions} 
\color{black} To solve the optimal control problem \eqref{ocp2}, we derive first-order conditions for optimality. These conditions are obtained by analyzing the behavior of $\mathcal{J}$ under infinitesimal perturbations of the control and, consequently, the state variable. Assuming $\mathcal{J}(u)$ is Fréchet differentiable (proof in Appendix \ref{differentiability}) with derivative denoted by $D\mathcal{J}(u)[\cdot]$, every optimal control signal $u^*\in L_2(0,T;L_\infty(\Omega))$ satisfies
\begin{equation}\label{PDE_optimality}
D\mathcal{J}(u^*)[h] = 0\qquad \forall h \in L_2(0,T;L_\infty(\Omega))\,.
\end{equation}

The condition \eqref{PDE_optimality} can be rewritten in a more convenient form by means of a system of equations for the triple $(u,\,\rho,\,\lambda)$, where $\lambda\in \mathcal{W}(0,T)$ denotes the adjoint state variable associated with the pair $(\rho,u)$
\begin{equation}\label{PMP_PDE}\begin{cases}
        \partial_t \rho = \beta^{-1}\Delta\rho + \nabla \cdot\bigg(\rho\big((\nabla W \ast \rho)+\nabla V + u\big)\bigg) &\text{(state)}\\[2ex]
        \color{black}-\partial_t \lambda = \beta^{-1}\Delta\lambda -\nabla\lambda\cdot \bigg(\nabla W*\rho + \nabla V + u\bigg) + \bigg(\nabla W\star(\nabla\lambda\,\rho)\bigg) + (\rho-\rho_\infty) &\text{(adjoint)}\\[2ex]
        \color{black}\rho(x,0) = \rho_0(x),\qquad\lambda(t,x) = \nabla_\rho\Psi(\rho(t,x)) &\text{(b.c.'s)}\\[2ex]
        -\nabla\lambda\rho + \gamma u = 0 &\text{(optimality)}
\end{cases}
\end{equation}
and where the operation $\star$ is defined as
\begin{equation*}
    (f\star g)(x) := \int\limits_{\Omega} f(y,x)g(y)dy\,.
\end{equation*}
This optimality system can be immediately derived from the first order condition \eqref{PDE_optimality} combined with the derivation of a computational expression for $D\mathcal{J}(u^*)[\cdot]$ \cite{PMP_PDE_troltzsch,robfleig,Borzi_secondorder,PMP_PDE_Aronna,PMP_PDE_Borzi}.

If an optimal control signal $u^*$ exists that solves \eqref{ocp2}, then the optimal triple $(\rho^*,\lambda^*,u^*)$, comprising the $u^*$-controlled state $\rho^*$, the corresponding adjoint $\lambda^*$, and $u^*$ itself, satisfy eq. \eqref{PDE_optimality}. For a comprehensive derivation and well-posedness analysis of optimality conditions for a similar control problem, we refer the reader to \cite{robfleig}.

\color{black}In the following, we discuss the construction of a numerical scheme for the approximation of such an optimal triple by following a discretize-then-optimize strategy. That is, instead of approximating the continuous optimality system \eqref{PDE_optimality}, we derive optimality conditions for a discretized version of \eqref{ocp2}, which has a similar forward-backward structure as in the continuous case. For this, following the method of lines, we consider the spatial approximation in $\mathcal{V}_L$ of the target steady state $\rho_\infty$, the distribution $\rho(t,\cdot)$ and the control signal $u(t,\cdot)$ at each fixed time $t$ as $L-$truncated Fourier series. A semi-discrete problem is then obtained by embedding the time dependency at the level of the Fourier coefficients
\begin{equation}\label{fp_c}
    \rho_\infty^L(x) = \sum\limits_{i=1}^L \tilde a_i \psi_i(x)\,,\qquad \rho^L(t,x) = \sum\limits_{i=1}^L a_i(t) \psi_i(x)\,,\qquad u^L(t,x) = \sum\limits_{i=1}^L\sum\limits_{j=1}^d u_{i,j}(t) \psi_i(x)\,.
\end{equation} 
\color{black}
The discretization maintains the assumption that the control signal $u^L(t,\cdot) \in L_\infty(\Omega)$ for all $t \in [0,T]$. This is ensured by the boundedness of our basis functions $|\psi_i(x)|{L\infty}\leq1$ $i=1,\dots,L$, and by $u^L$ being a finite sum of such bounded functions, hence remaining in $L_\infty$.
\color{black}
The associated L-dimensional approximation of \eqref{ocp} reads
\begin{equation}\label{discrete_dyn}\color{black}
    M\,\dot{\mathbf{a}}(t) = -\big(\beta^{-1} A+C\big)\,\mathbf{a}(t) - b(\mathbf{a}(t)) - d(\mathbf{a}(t),\mathbf{u}(t))\,,
\end{equation} 
where the operators $\,A,\,b(\cdot)$ are as in (\ref{operators}-\ref{bil}), $\mathbf{u} \in \R^{L\times d}$, and the mass matrix $M$ has entries
\begin{equation*}
    M_{n,m} = \int\limits_{\Omega} \psi_n(x)\psi_m(x) dx\,,
\end{equation*}
whilst the bilinear control operator reads
\begin{equation*}
    d(\mathbf{a},\mathbf{u}) = \sum\limits_{i,k=1}^L\sum\limits_{j=1}^d a_i(t)\,u_{k,j}(t) \int\limits_{\Omega}\psi_i\frac{\partial\psi_m}{\partial x_j}\psi_k \,dx,\qquad m=1,\dots,L\,.
\end{equation*}

Under this discretization, denoting the $j-$th column of the control coefficients as $\mathbf{u}_j\in \R^L$, the finite-dimensional approximation of the optimal control problem reads
\begin{equation}\label{ocpd}
\underset{\mathbf{u}(\cdot)\in\mathcal{V}^d_L}{\min}\cJ_L(\mathbf{u}):=\dfrac{1}{2}\int\limits_0^T (\mathbf{a}(t)-\tilde{\mathbf{a}})^\top M\, (\mathbf{a}(t)-\tilde{\mathbf{a}}) + \gamma \sum\limits_{j=1}^d (\mathbf{u}_{j}(t)^\top M\, \mathbf{u}_{j}(t)) \,dt + P(\mathbf{a}(T))\,,    
\end{equation}
subject to \eqref{discrete_dyn}, where final cost is given by
\begin{equation*}
    P(\mathbf{a}(T)) = \eta\,(\mathbf{a}(T)-\tilde{\mathbf{a}})^\top M\, (\mathbf{a}(T)-\tilde{\mathbf{a}})\,.
\end{equation*}

For the sake of simplicity, we will drop the time dependence and denote $\mathbf{a}=\mathbf{a}(t),\;\mathbf{u} = \mathbf{u}(t)$. Writing the system dynamics in \eqref{discrete_dyn} in control-affine form as
\begin{equation*}
    \dot{\mathbf{a}} = f(\mathbf{a}) + \sum\limits_{j=1}^d g_j(\mathbf{a})\,\mathbf{u}_j
\end{equation*}
\begin{equation*}
    f(\mathbf{a}) := M^{-1}\bigg(-(\beta^{-1}A+C)\,\mathbf{a} - b(\mathbf{a})\bigg)\,,\quad g_j(\mathbf{a}):= M^{-1}\bigg(-D_j(\mathbf{a})\bigg)\,, \quad \color{black}g(\mathbf{a}) := \big(g_1(\mathbf{a}),\dots,g_d(\mathbf{a})\big),\color{black}
\end{equation*}
\begin{equation*}
\big[D_j\big]_{k,m}(\mathbf{a}) := \sum\limits_{n=1}^L a_n\int\limits_\Omega \psi_n \frac{\partial\psi_m}{\partial x_j} \psi_k\, dx\,, \qquad k,m = 1,\dots,L,
\end{equation*} 
we define the Hamiltonian associated to the optimal control problem
\begin{equation*}
    \mathcal{H}(\mathbf{a},\mathbf{u},\mathbf{p}):=\mathbf{p}^\top\bigg(f(\mathbf{a}) + g(\mathbf{a})\,\mathbf{u}\bigg)+l(\mathbf{a},\mathbf{u})\,, 
\end{equation*} 
with adjoint variable $\mathbf{p}\in\R^L$ and  running cost 
\begin{equation*}
    l(\mathbf{a},\mathbf{u}):=(\mathbf{a}-\tilde{\mathbf{a}})^\top M\, (\mathbf{a}-\tilde{\mathbf{a}})+\gamma \,\sum\limits_{j=1}^d (\mathbf{u}_{j}^\top M\, \mathbf{u}_{j})\,.
\end{equation*}
Analogously as for the continuous problem, first-order necessary optimality conditions for the semi-discrete problem \eqref{ocpd} are then expressed as
\begin{equation*}
    \dot{\mathbf{a}} = \frac{\partial \mathcal{H}}{\partial \mathbf{p}}(\mathbf{a}, \mathbf{u},\mathbf{p})\,,\quad
    -\dot{\mathbf{p}} = \frac{\partial \mathcal{H}}{\partial \mathbf{a}}(\mathbf{a}, \mathbf{u},\mathbf{p})\,,\quad
    \frac{\partial \mathcal{H}}{\partial \mathbf{u}}(\mathbf{a}, \mathbf{u},\mathbf{p})=0\,,
\end{equation*}
closed with initial and terminal boundary conditions. This two-point boundary value problem reads
\begin{align}
        \dot{\mathbf{a}} &=\; M^{-1}\bigg(-(\beta^{-1} A+C)\,\mathbf{a} - b(\mathbf{a}) - d(\mathbf{a},\mathbf{u})\bigg)&\quad\text{(state)}\label{state}\\
        -\dot{\mathbf{p}} &=\; \mathbf{p}^\top M^{-1}\nabla_{\mathbf{a}} \bigg(-(\beta^{-1} A+C)\,\mathbf{a} - b(\mathbf{a}) - d(\mathbf{a},\mathbf{u})\bigg)+ M\,(\mathbf{a}-\tilde{\mathbf{a}})&\quad\text{(adjoint)}\label{adjoint}\\[2ex]
        \mathbf{a}(0) &=\; \hat{\mathbf{a}}\,,\quad\mathbf{p}(T)=\nabla_\mathbf{a} P(\mathbf{a}(T)) = 2\eta\,(\mathbf{a}(T)-\tilde{\mathbf{a}})^\top M&\quad\text{(b.c.'s)}\label{bcs}\\[2ex]
        0&=\gamma\, \mathbf{u}_j^\top M - \mathbf{p}^{\top} g_j(\mathbf{a}), \qquad j = 1,\dots,d &\quad\text{(optimality)}\label{TPBVP}      
\end{align}
The system is addressed by combining the forward-backward integration of  the state \eqref{state} and adjoint \eqref{adjoint} equations, with a gradient descent algorithm for the update of the control signal. Indeed, the optimality condition \eqref{TPBVP} corresponds to $\nabla_{\mathbf{u}}\cJ_L=0$. This iterative procedure is shown in Algorithm \ref{alg:PMP}.

\begin{algorithm}\caption{Reduced gradient iteration}\label{alg:PMP}
\KwData{initial condition $\hat{\mathbf{a}}$; initial guess for the control signal  $\mathbf{u}^1$; time horizon $T$; \\$\qquad\;\;$ tolerance $\varepsilon\ll1$; step size $\delta > 0$}
\KwResult{$(\mathbf{a}, \mathbf{u},\mathbf{p})$ approximate solution of optimality system (\ref{state} - \ref{TPBVP})}
$n=0$, converged = \texttt{false}\\
\While{not converged}{
$n = n+1$\\
$\mathbf{u}\;\,\,\gets\mathbf{u}^n$\\
$\mathbf{a}^n \gets$ integration of \eqref{state} from $\hat{\mathbf{a}}$ forward in time \\
$\mathbf{p}^n(T)\gets$ \eqref{bcs}\\
$\mathbf{p}^n \gets$ integration of \eqref{adjoint} from $\mathbf{p}^n(T)$ backward in time \\
$\mathbf{u}^{n+1} \gets \mathbf{u}^n - \delta\bigg(\gamma\, \mathbf{u}^{n\top}M - g_j(\mathbf{a}^n) \bigg)$\Comment*[r]{gradient step}
\emph{converged} \textbf{if} $\|\mathbf{u}^{n+1}-\mathbf{u}^{n}\|_2^2 < \varepsilon$
}
\end{algorithm}

\color{black}\begin{remark} (On the consistency of  optimality conditions). The continuous optimal control problem \eqref{ocp} and its discretized counterpart \eqref{ocpd} yield different optimality systems. 
Establishing consistency between these two approaches is essential for ensuring the reliability and accuracy of our numerical solutions. We provide a detailed proof of this consistency in Appendix \ref{consistency_appendix}.\end{remark}
\color{black} 

\subsection{Nonlinear Model Predictive Control}

Our original goal is the design of a control law for stabilizing the McKean-Vlasov PDE towards an unstable steady state. In general, optimal asymptotic stabilization is addressed by means of infinite horizon optimal control, that is, taking $T=+\infty$. This formulation is not approachable with the proposed methodology as it would translate into an ill-posed two-point boundary value problem with an asymptotic terminal condition as $t\to\infty$. 

On the other hand, finite horizon $T<+\infty$  optimal control problems are generally not suitable for stabilizing systems around unstable stationary configurations. Small perturbations or deviations from these points tend to derail the system from the desired convergence and in our case, trajectories will be attracted towards stable steady states. Even if the open-loop optimal control successfully steers the system towards an unstable steady state $\rho_\infty$ within the finite horizon, the system is likely to diverge from that beyond $T$ due to disturbances. In this context, a fundamental aspect is the use of a feedback control law which can continuously steer the dynamics towards the desired configuration based on measurements of the current state of the system, hence robust to perturbations. 

Nonlinear Model Predictive Control (MPC) is a technique for synthesizing stabilizing feedback control maps from subsequent solves of the optimality system \eqref{state}-\eqref{TPBVP} over a moving time frame of reduced horizon $T<+\infty$. 
The idea is that state dependency of the open loop solution only resides in the first entry in the control signal, which can be registered as feedback w.r.t. the initial condition. This procedure is presented in Algorithm \ref{alg:MPC}. The receding horizon nature of MPC allows for subsequent feedback and correction of deviations from the unstable target configuration.

\begin{algorithm}\caption{MPC}\label{alg:MPC}
\KwData{$\hat{\mathbf{a}}$ initial condition; $T$ reduced time horizon; $N\_steps$ number of time steps in $[0,T]$}
\KwResult{Feedback controlled trajectory $\mathbf{a}_h$, $h=1,\dots,N\_steps$}
$dt = 1/(N\_steps+1)$\\
$\mathbf{a}_0 \gets \hat{\mathbf{a}}$\\
\For{$h = 1,\dots,N\_steps$}{
$(\mathbf{a}(t), \mathbf{u}(t),\mathbf{p}(t))\gets$ Algorithm \ref{alg:PMP} with initial condition $\mathbf{a}_{h-1}$ and horizon $T$\\
$\mathbf{u}_h \gets \mathbf{u}(t=0)$\Comment*[r]{feedback control}
$\mathbf{a}_{h} \gets $ integration of \eqref{state} from $\mathbf{a}_{h-1}$ for $t\in[0,dt]$ with control $\mathbf{u}_h$
}
\end{algorithm}

\section{Numerical Tests}
\label{sec:numerics}
In this section we assess the proposed deflation and control methodology on a test cases related to opinion dynamics. In this context, the uncontrolled dynamics converge to a single cluster (or consensus), and we wish to prevent such behaviour and steering the system towards an unstable multi-cluster or uniform distribution. 

{\color{black}The PDE model parameters, such as noise and interaction strengths, are selected to ensure the system operates below the transition point where multiple solutions (both stable and unstable) to the integral equation emerge.  For a comprehensive theoretical analysis of these models, including the examination of critical temperatures and interaction strengths, readers are referred to~\cite{GPY2017, CGPS2020} for the Hegselmann-Krause model and in~\cite{GPS_2024} for the 2d Von Mises model.}

In both the following numerical tests, at a fixed time, the control signal and density are discretized in space as in \eqref{fp_c} as linear combination of spectral basis functions, here given by the first $L$ Fourier modes. For the time stepping we will use a $4$th order Runge Kutta scheme.

\subsection{Hegselmann-Krause opinion formation model}
We consider the controlled McKean–Vlasov equation with $V\equiv0$ and Hegselmann-Krause interaction potential \begin{equation*}
    \nabla W(x,y):=\begin{cases}
        x-y &\text{if } |x-y|\leq r\\
        0 &\text{otherwise}
    \end{cases}
\end{equation*}
where $r$ encodes the radius of interactions. The jump discontinuities in this potential are treated by shifting either the base or the end of the jump by a small factor of $5e-3$. 

\begin{figure}[h]
\centering
        \includegraphics[width=0.3\textwidth]{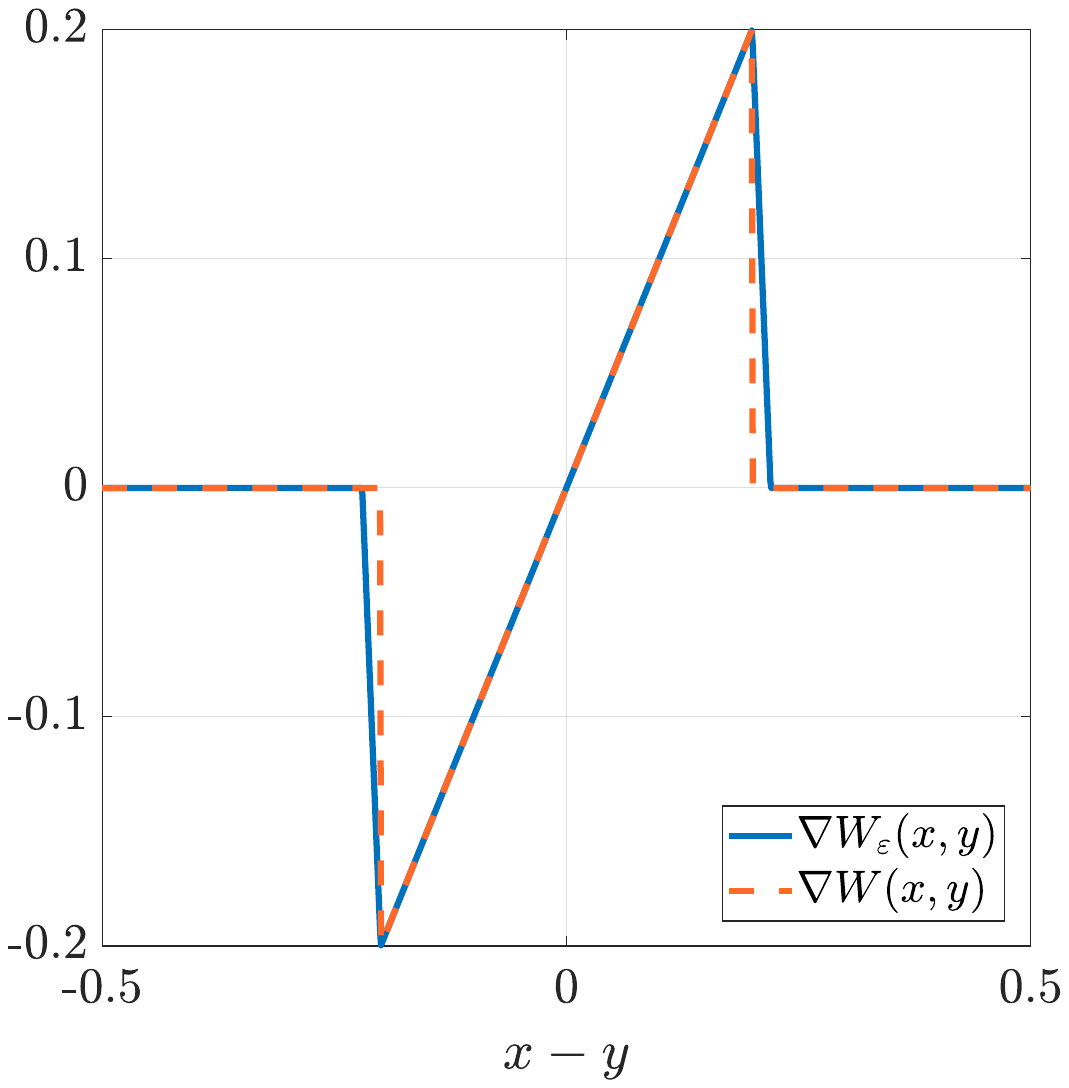} \includegraphics[width=0.3\textwidth]{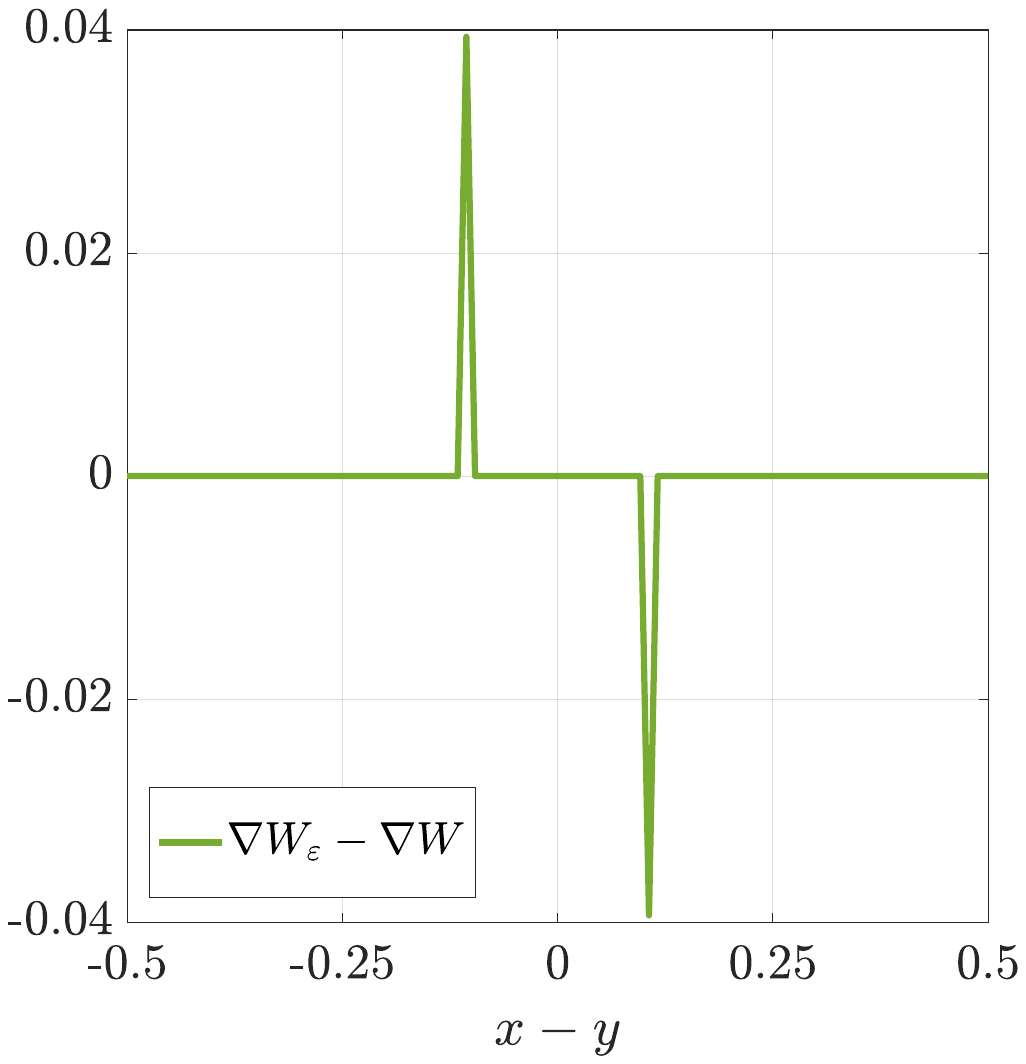}
        \caption{\color{black}The Hegselmann-Krause interaction potential and its regularization (left). Difference between original and regularized potential (right).}
    \end{figure}  
    
In the one-dimensional case with $\Omega = [0,1]$ and periodic boundary conditions, the Mckean-Vlasov equation reads
\begin{equation*}
    \partial_t \rho = \beta^{-1} \partial^2_{x}\rho +  \partial_x\bigg[\rho\bigg(\int_{\Omega} (x-y)\rho(y)\mathbb{I}_{|x-y|\leq r}(y)\,dy\bigg)\bigg]\,,
\end{equation*}
where $\mathbb{I}_{|z|\leq r}$ denotes the indicator function of the ball of radius $R$ centered at the origin. 

\begin{figure}[!h]
    \centering
    \includegraphics[width = 0.32\textwidth]{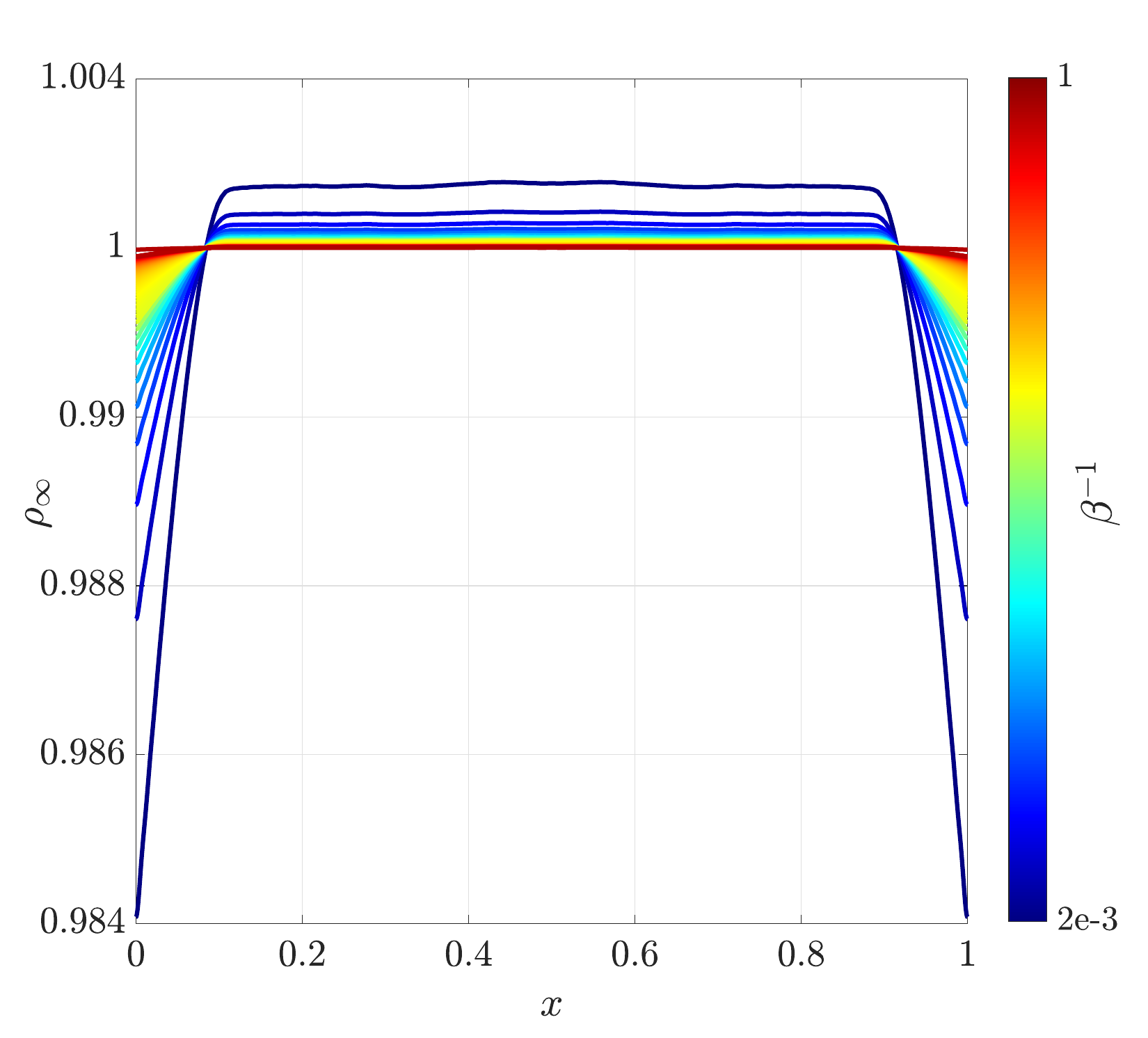}\hfill
    \includegraphics[width = 0.3\textwidth]{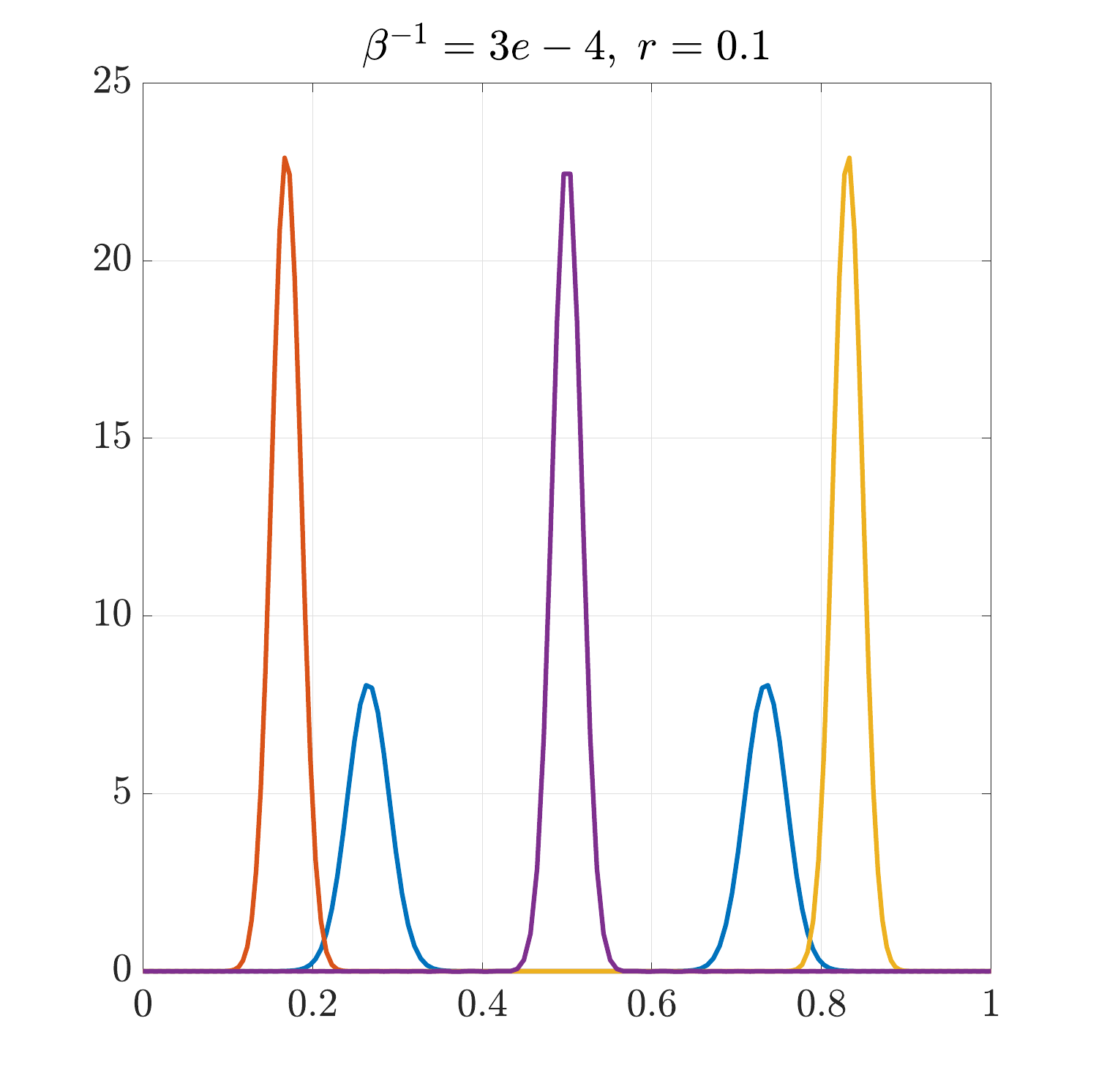}\hfill
    \includegraphics[width = 0.3\textwidth]{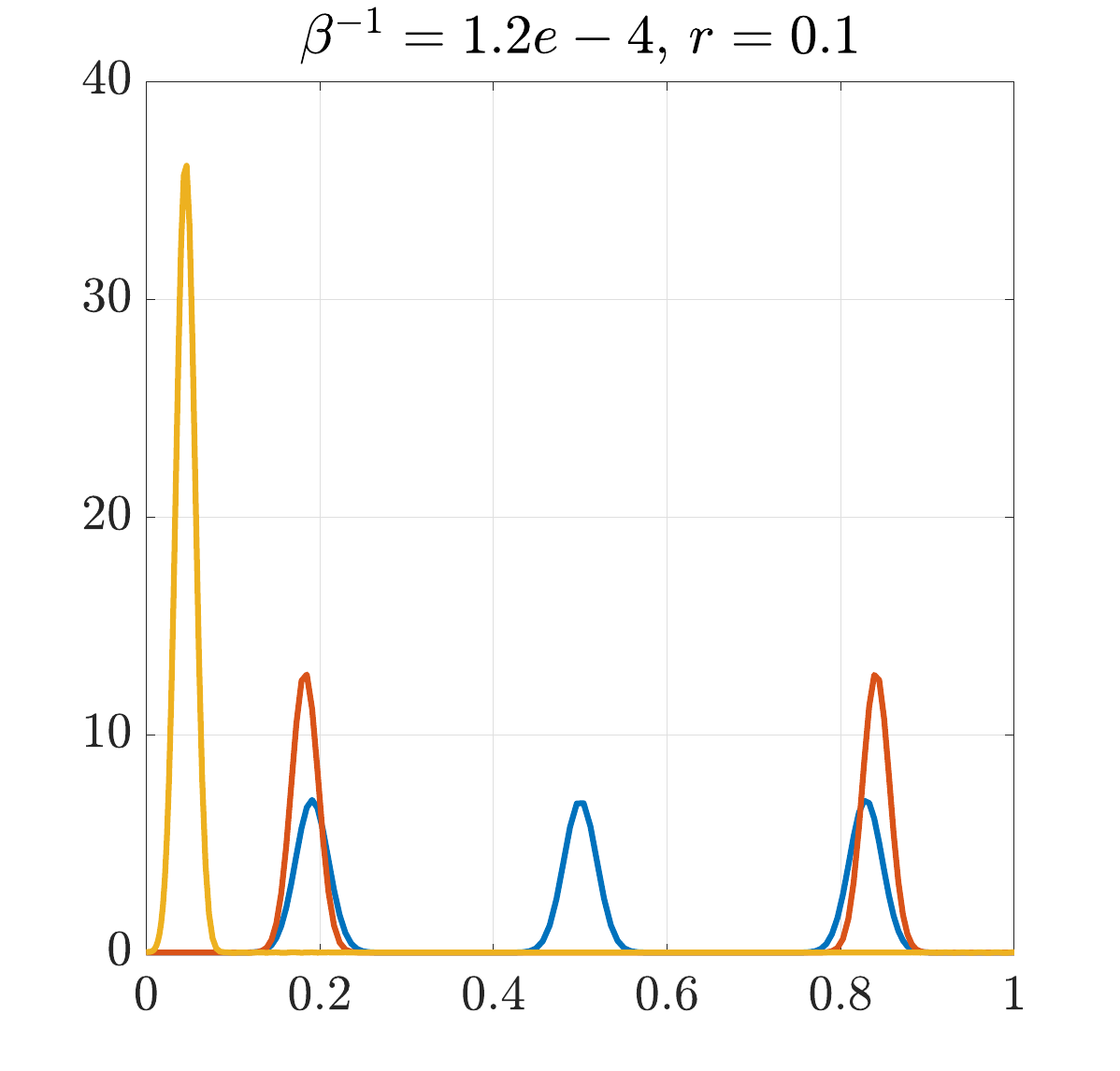}
    \caption{In high noise regimes ($\beta^{-1}>\beta_c^{-1}\approx 1e-3$), a single stable steady state converges towards the disordered uniform state (left). For low noise regimes (centre-right), opinions start condensing into clusters. In both parameter regimes, the stable configuration is consensus. However, deflated Newton's method identifies some unstable steady states displaying concentration around multiple positions. All solutions can be found translated along the domain, but only a selection of those is displayed for clarity of presentation.}
    \label{Fig:ss_opinions}
\end{figure}

The emergent behaviour in this model depends on both the noise amplitude $\beta^{-1}$ and the range of interaction $r$. In Figure \ref{Fig:ss_opinions} we fix $r=0.1$ and we analyze the asymptotic behaviour while changing $\beta^{-1}$. The figures have been obtained with $l=50$ Fourier modes, $200$ Gauss points, and deflation parameters $p=3$, $\xi=1$. 
In high noise regimes (similarly to what happens for long interaction ranges), we identify a unique stationary configuration, converging to the uniform distribution as $\beta^{-1}\to 1$.
For lower inverse temperature $\beta^{-1}\leq\beta^{-1}_c$, agents' opinions tend to cluster into distinct groups~\cite{GPY2017}. Once the opinion clusters are formed, they remain relatively stable for a long period of time. This is a manifestation of the dynamical metastability phenomenon that characterizes the dynamics of systems exhibiting discontinuous phase transitions, and it will be studied in detail in forthcoming work. Eventually, the opinion clusters merge, leading to the stable stationary consensus configuration. 

\begin{figure}[h]
    \centering
    \includegraphics[width = 0.32\textwidth]{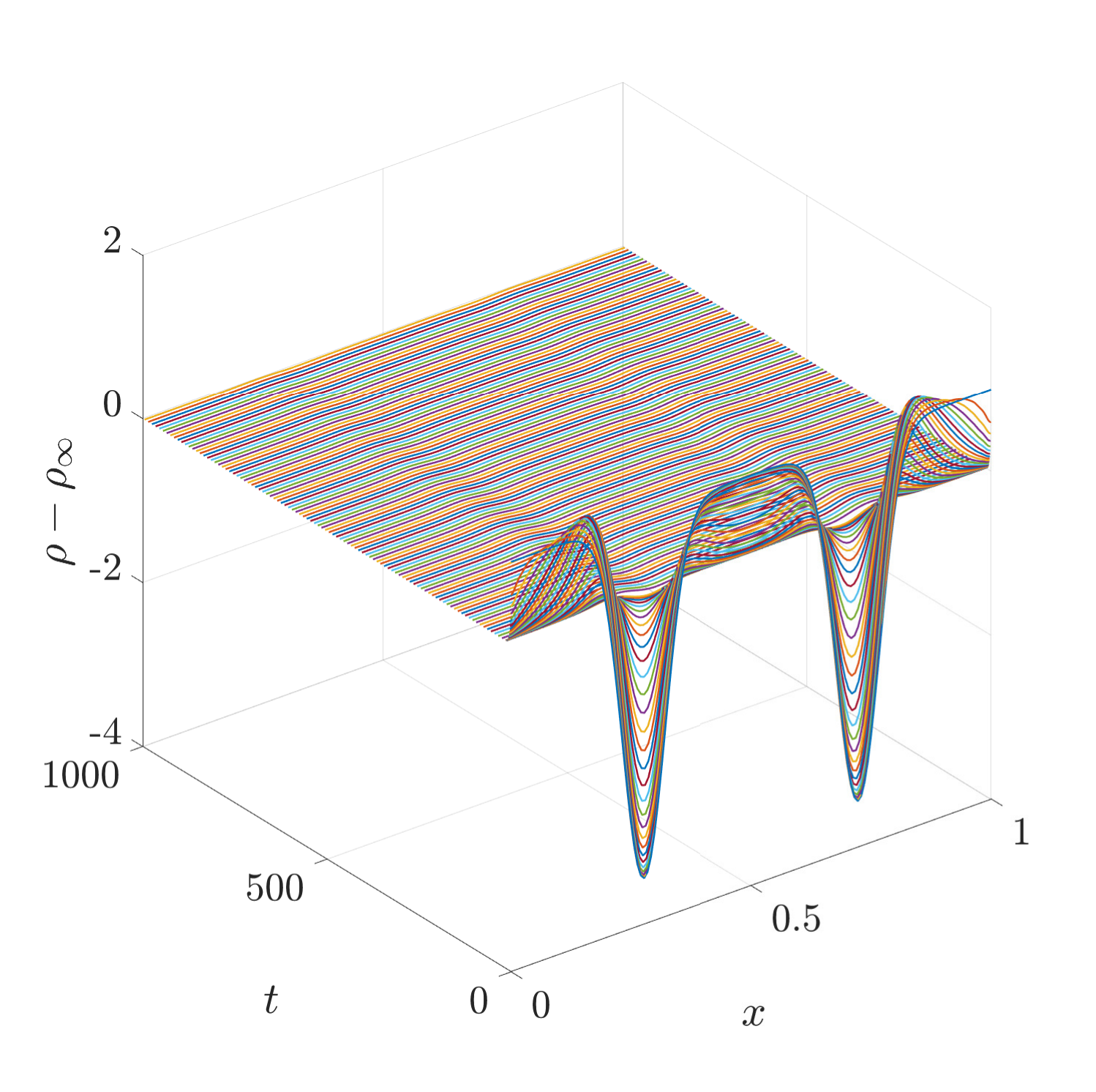}
    \includegraphics[width = 0.32\textwidth]{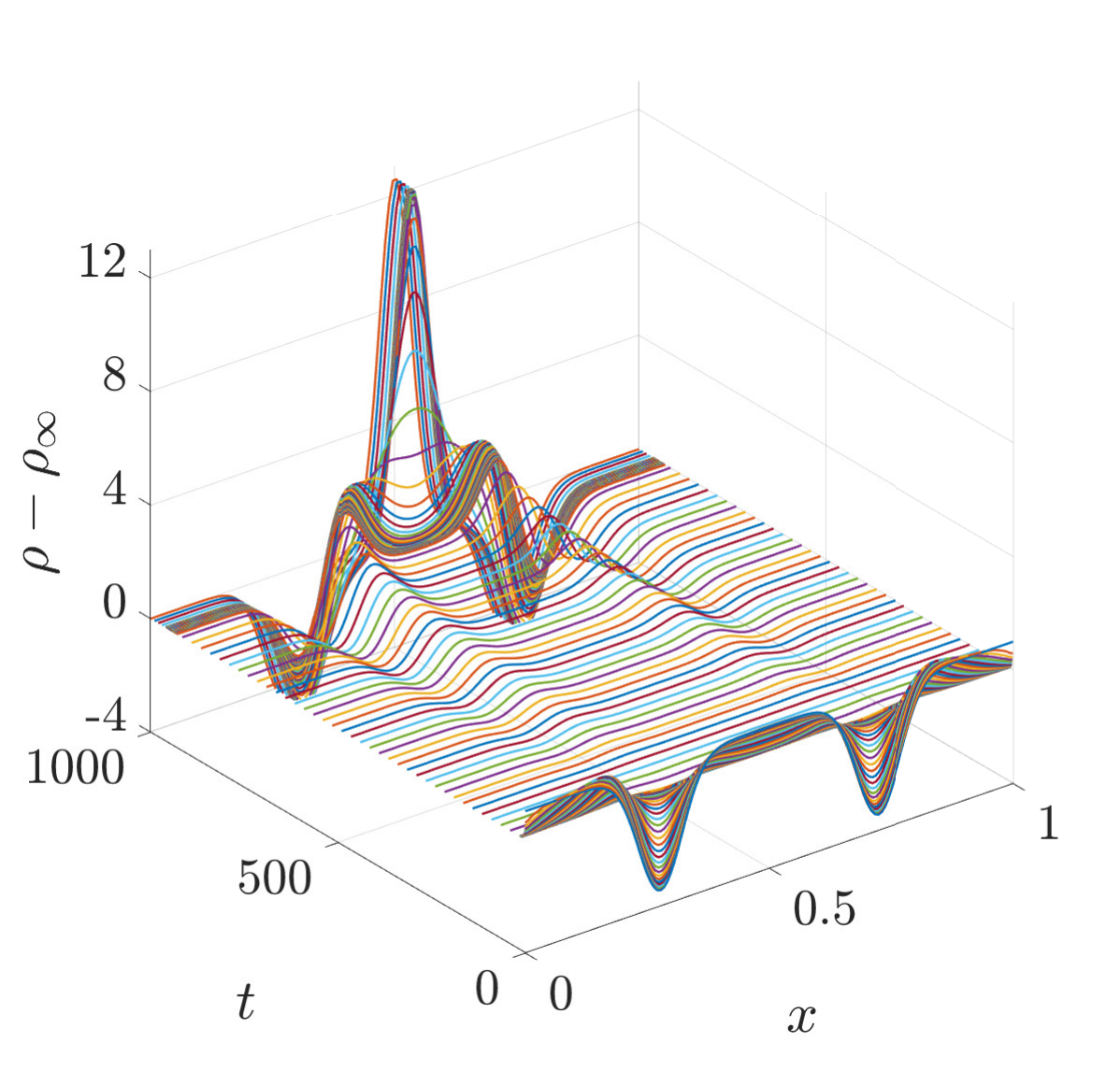}
    \includegraphics[width = 0.3\textwidth]{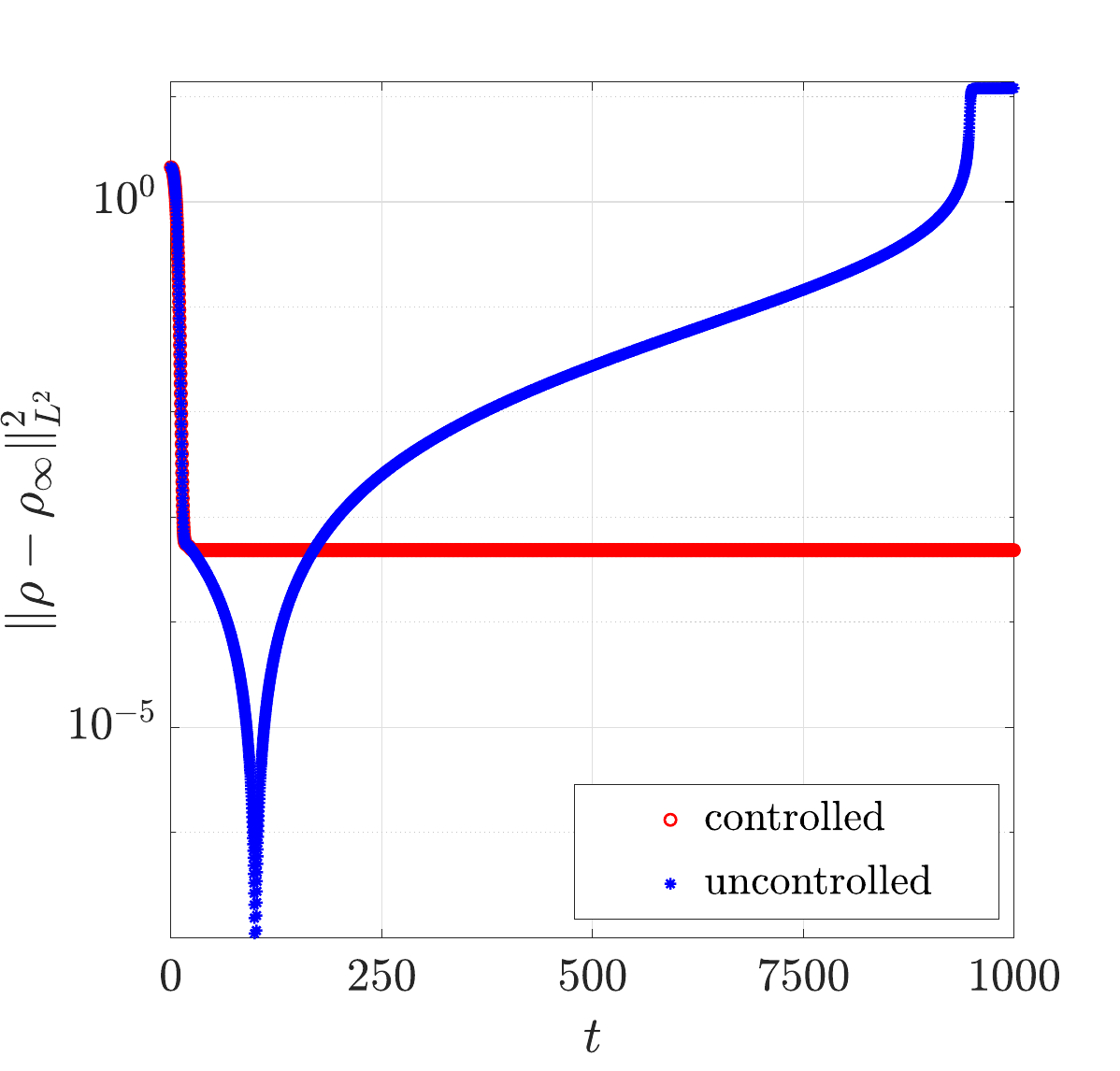}
    \caption{We address the optimal control problem for stabilizing the system in (Fig. \ref{Fig:ss_opinions} - centre) at the unstable configuration in blue via MPC.  (Left) evolution of the discrepancy w.r.t. $\rho_\infty$ in the controlled case, departing from a uniform initial condition. 
    (Centre) the free dynamics diverge from the target state, towards the stable consensus configuration.
    (Right) discrepancy in $L^2$ norm between density $\rho(t,x)$ and target steady state $\rho_\infty(x)$.}
    \label{Fig:ss_opinions_control}
\end{figure}

Having identified the multi-cluster distributions as metastable, we solve the optimal control problem \eqref{ocp} for stabilizing the system at this configuration while preventing consensus formation. We consider   $\beta^{-1}=3e-4$, with the stabilization target distribution $\rho_\infty$ being the blue line in Figure \ref{Fig:ss_opinions} (centre). 
We address the stabilization problem with control energy penalization parameter $\gamma=0.01$, and terminal cost weight $\eta = 100$, over a total time interval $[0, T ]$, $T = 1e3$ seconds. {\color{black}The parameters for the gradient descend are a step size $\delta=10^{-6}$ and convergence tolerance $\varepsilon=0.01$.} In a receding horizon approach, the control is computed by subsequent solves of \eqref{TPBVP} over a rolling window of reduced time horizon ($100$s). In Figure \ref{Fig:ss_opinions_control}, we compare the controlled and uncontrolled evolution of $\rho(t,x)$, where the MPC action successfully stabilizes the distribution at the desired state, preventing the formation of the concentration profile appearing in the uncontrolled dynamics.


\subsection{2d Von Mises interaction potential}
\begin{figure}[!h]
\centering
\begin{subfigure}[h]{0.33\linewidth}
    \centering
    \includegraphics[width = 1\textwidth]{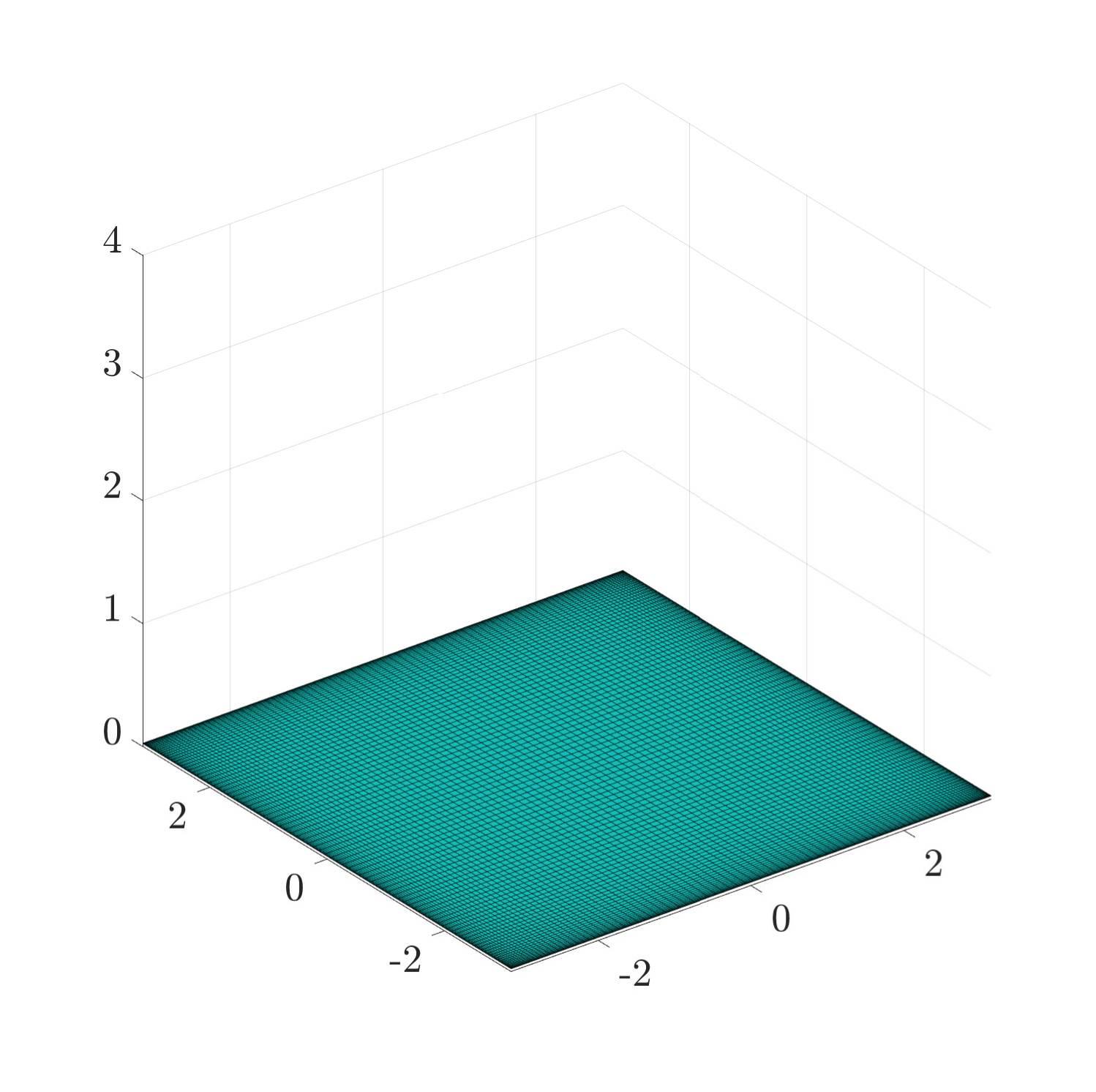}
    \caption{$\beta^{-1} > 0.3702$}
\end{subfigure}
\begin{subfigure}[h]{0.33\linewidth}
    \centering
    \includegraphics[width = 1\textwidth]{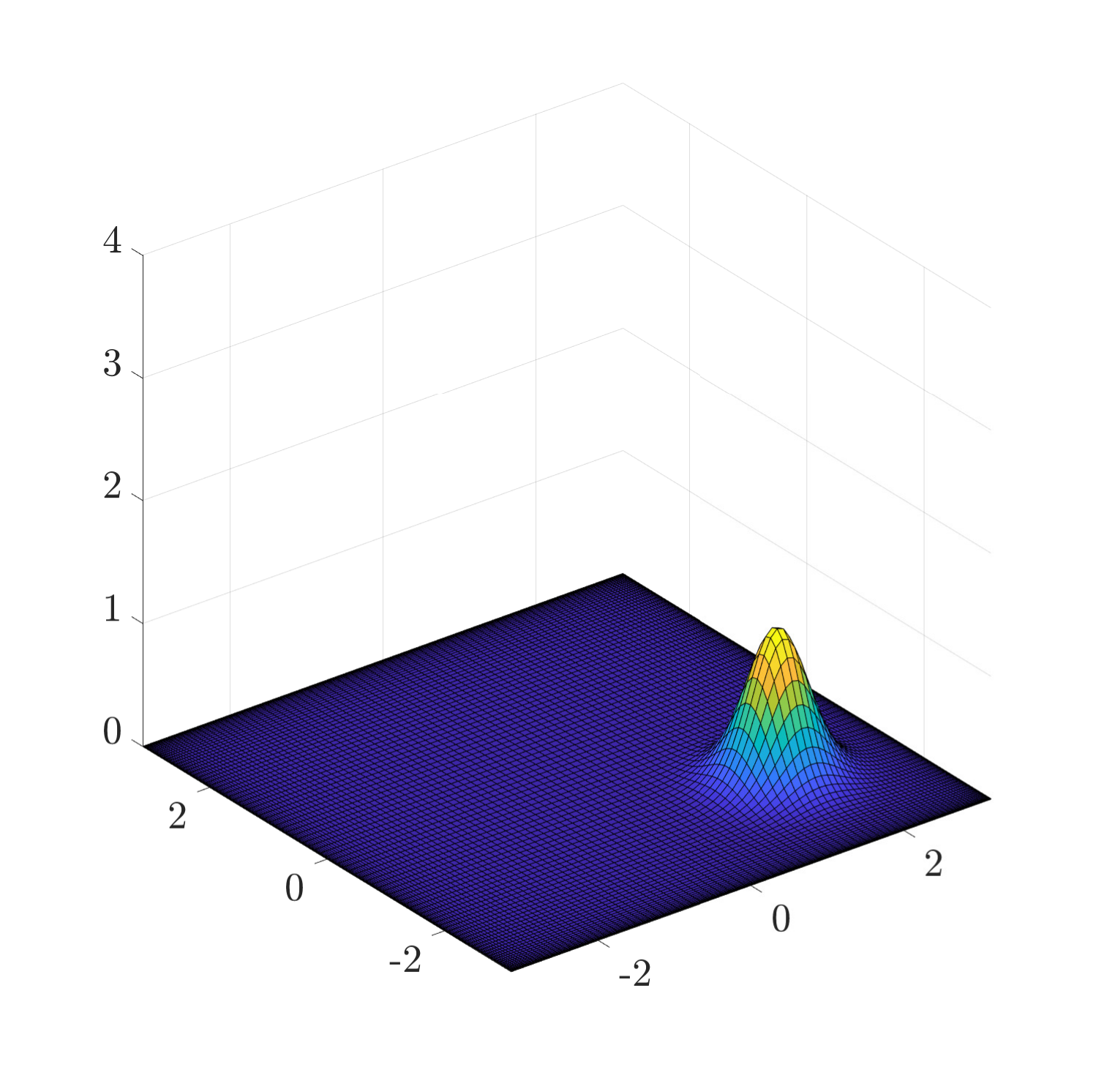}
    \caption{$\beta^{-1} = 0.3701$}
\end{subfigure}
\begin{subfigure}[h]{0.32\linewidth}
    \centering
    \includegraphics[width = 1\textwidth]{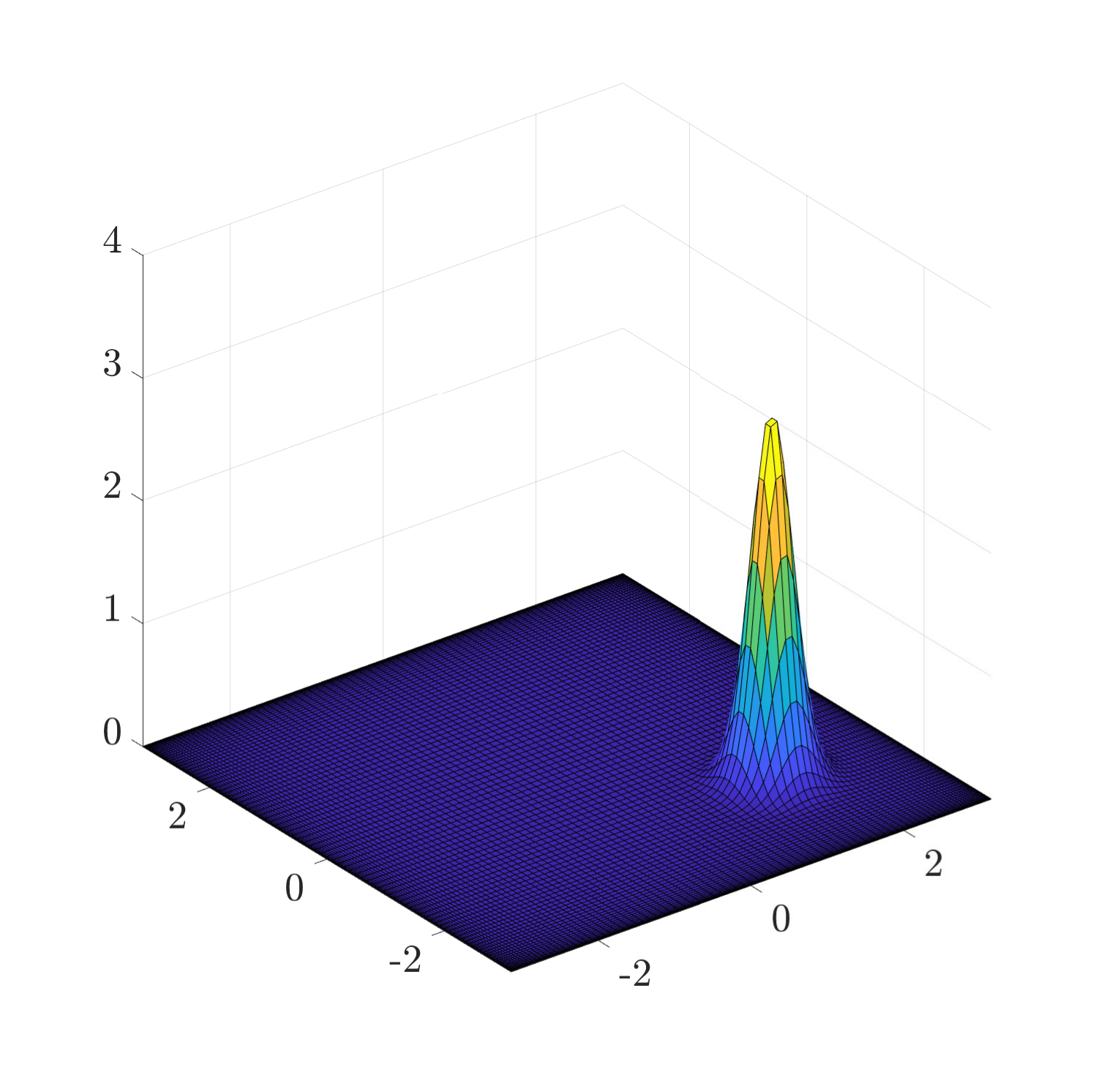}
    \caption{$\beta^{-1}= 0.2$}
\end{subfigure}
    \caption{Steady states identified with the deflated Newton's method. For $\beta\geq0.3702$, the iterative scheme converges into a uniform distribution, and the root finding method does not identify other steady states. For smaller values of $\beta^{-1}$, we have convergence towards a droplet configuration, whilst the uniform distribution remains an unstable steady state.}
    \label{fig:deflation_VonMises}
\end{figure}

In this numerical test we apply the proposed methodology to a two dimensional problem, with $\Omega=[-\pi,\pi]\times[-\pi,\pi]$ and Von Mises interaction potential
\begin{equation*}\label{e:von_mises}
     W(x,y) = -I_0(k)^{-2}exp\bigg\{\theta\big(cos(x)+cos(y)\big)\bigg\},
\end{equation*}
without any confining force $V(x)\equiv0$. We use this a "periodized" version of the two-dimensional attractive Gaussian interaction potential that was used in, e.g.~\cite{Martzel_Aslangul_2001, GRV-E_2022}.\footnote{Alternatively, we could have also considered the wrapped Gaussian. For small variances, the case of interest for us, the wrapped Gaussian, the Gaussian potential and the Von Mises potential give identical numerical results. On the other hand, the rigorous theory put forward in~\cite{CGPS2020} upon which we rely in order to prove the existence of a discontinuous phase transition for the bounded confidence model, applies only to periodic interaction potentials.} For sufficiently small interaction parameter, relative to the size of the domain, this system displays a discontinuous phase transition between a disordered phase, where at high temperatures $\beta^{-1}\geq\beta_c$ the stable steady state is the uniform distribution, and a cooperative one, where the uniform configuration is still a stationary state (now unstable), together with a stable, "droplet" concentration profile~\cite{CGPS2020}. In this numerical test we fix $\theta=1$. 
In the deflation routine, as initial guess $\mathbf{a}_0$, we consider the one resulting from the iterative scheme presented in \cite{GRV-E_2022}. As expected, for large values of $\beta^{-1}$, we have a single (uniform) solution, whilst the solutions become two when the inverse temperature is low enough. Figure \ref{fig:deflation_VonMises} shows the results obtained for $\kappa=1$, with $l=10$ Fourier modes and $200$ Gauss points. 

\begin{figure}[!h]
    \centering
    \begin{subfigure}{0.32\textwidth}
        \includegraphics[width = 1\textwidth]{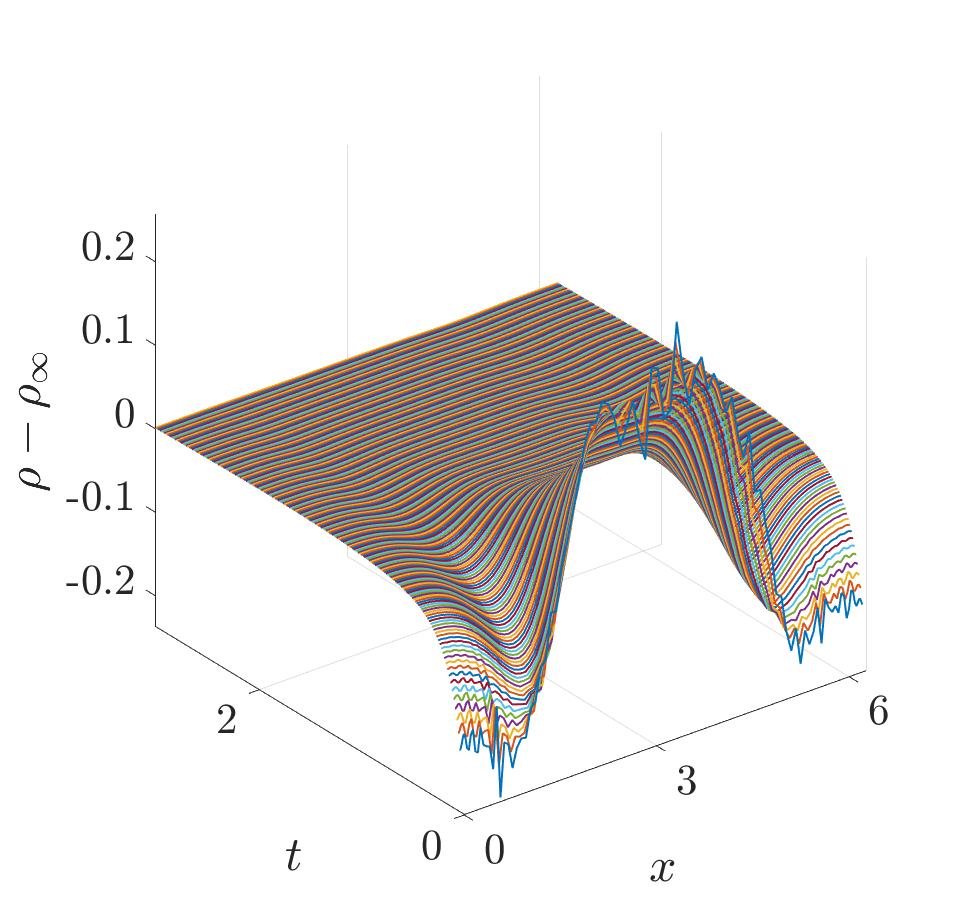}\vspace{-1cm}\caption{ }
    \end{subfigure}
    \hfill
    \begin{subfigure}{0.32\textwidth}
        \includegraphics[width = 1\textwidth]{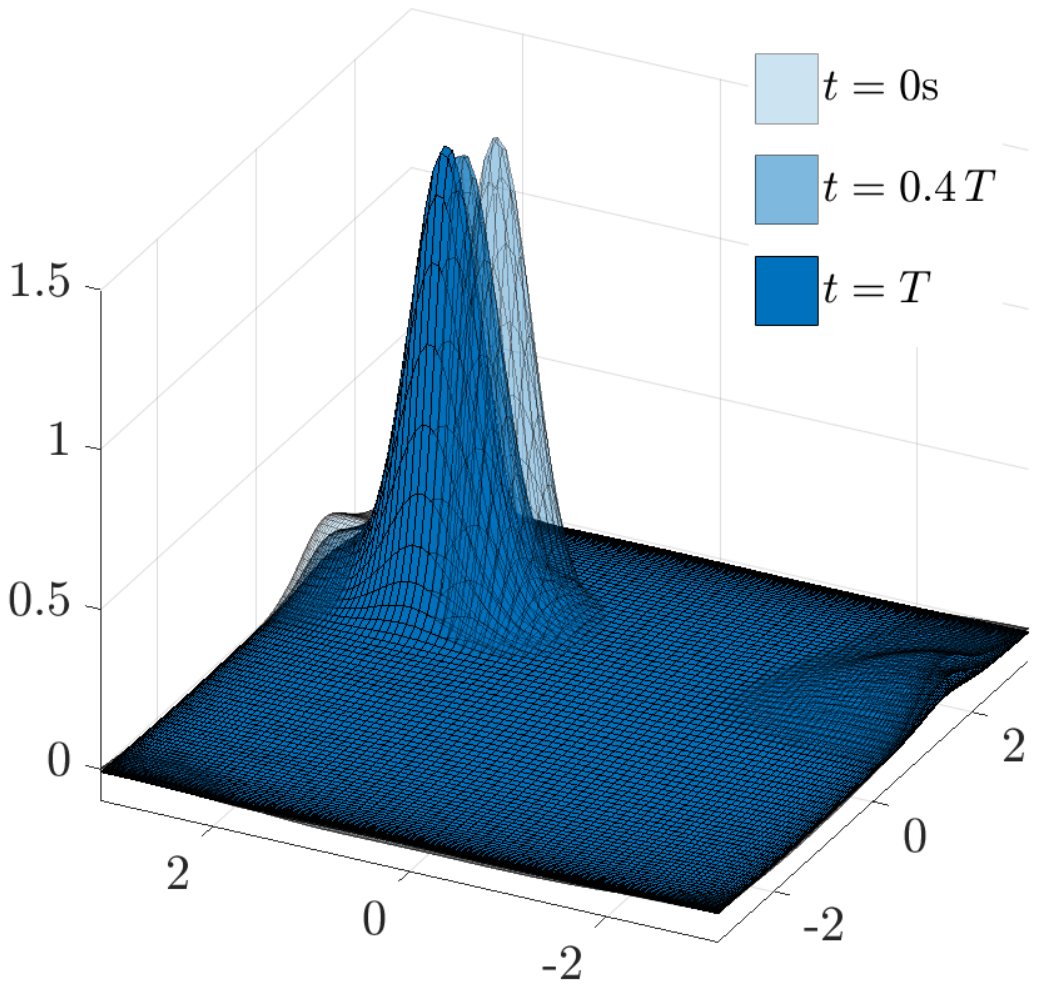}\vspace{-1cm}\caption{}
    \end{subfigure}\hfill
    \begin{subfigure}{0.32\textwidth}
        \includegraphics[width = 0.87\textwidth]{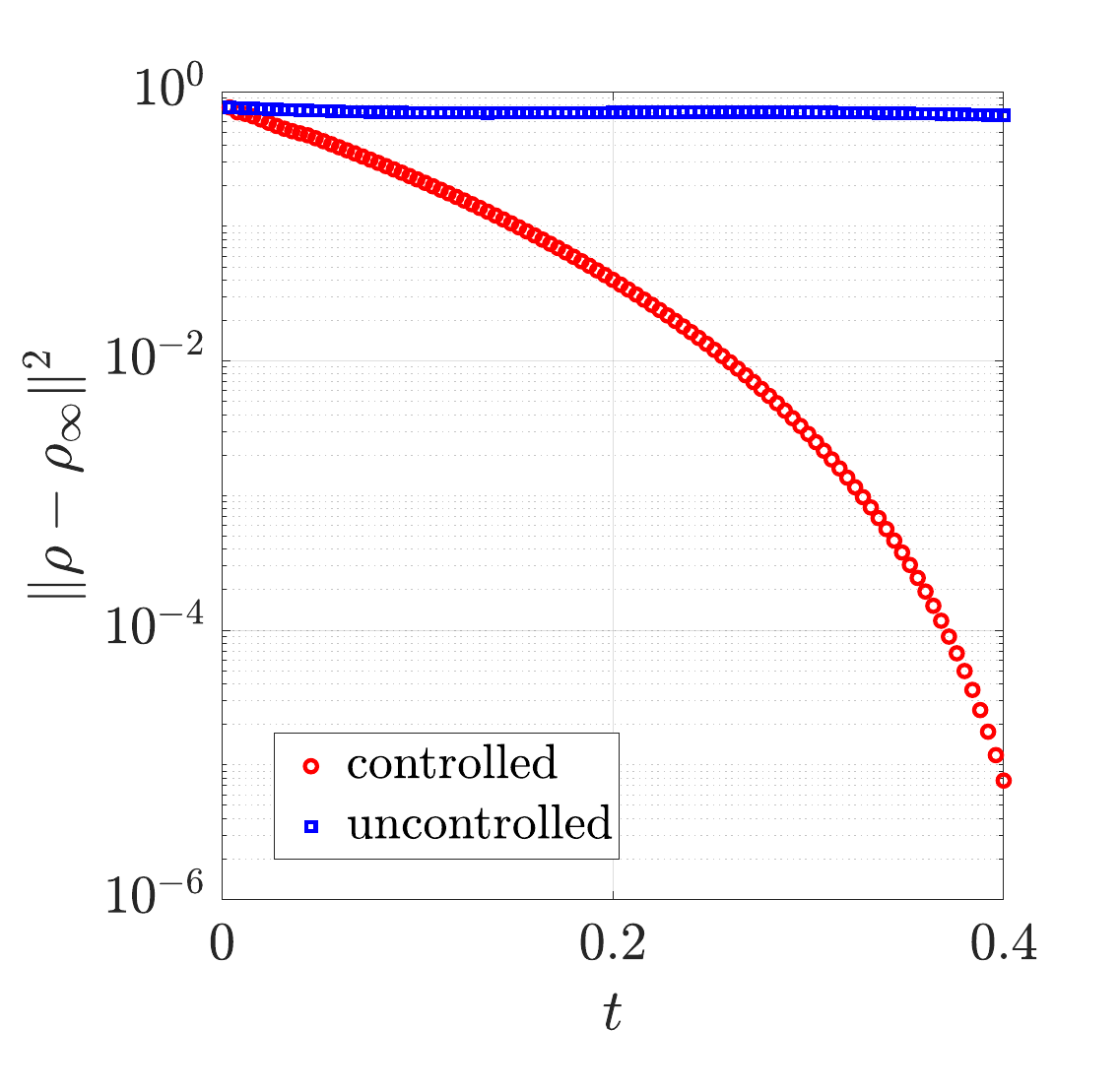}\caption{ }
    \end{subfigure}
    \caption{(a) evolution of the controlled dynamics towards $\rho_\infty$ at three sample times. 
    (b) the free dynamics diverge from the target state, forming the stable droplet configuration.
    (c) discrepancy in $L_2$ norm between density $\rho(t,x)$ and target steady state $\rho_\infty(x)$.}
    \label{fig:control_VonMises}
\end{figure}

We consider an optimal control problem for stabilizing the system towards the uniform distribution in a parameter regime where this is an unstable steady state ($\beta^{-1} = 0.3701$). Following a receding horizon strategy, we solve the problem in $t\in[0,T]$, $T =1s$ with a rolling window of reduced $0.1$s length, for control penalization parameter $\gamma=0.001$ and terminal cost weight $\eta=1000$. {\color{black}The parameters for the gradient descend are  $\delta=10^{-3}$ and $\varepsilon=0.1$.} In Figure \ref{fig:control_VonMises} we compare the evolution of $\rho(t,x)$ with and without the control action.  In the controlled case, the distribution of particles is successfully steered to the target uniform distribution, whilst the free dynamics move towards the stable concentration profile.

\section{Conclusions}\label{sec:conclusions}
We have presented a comprehensive methodology for computation and control of unstable steady states for McKean-Vlasov PDEs. The proposed pipeline consists of a spectral discretization of the dynamics, the use of a deflation algorithm for identification o all possible steady states, and a nonlinear model predictive control law for asymptotic stabilization towards the desired unstable state.  While unstable steady states do not naturally arise from studying the evolutionary system, there are applications where it might be useful to artificially force the system around an unstable configuration. For example, in opinion dynamics, one may wish to break down clustering in a highly polarised society. In this context, we have demonstrated the effectiveness of our method in multiagent models for opinion dynamics both in one and two dimensions. 

Even though in this paper we applied our methodology to a relatively simple class of examples, namely weakly interacting diffusions on the torus, for which a well developed theory of phase transitions exists, our approach opens up new possibilities for understanding and controlling the collective behavior of complex interacting particle systems, with potential impact in various scientific and engineering domains. We are particularly interested in multiagent systems in the social sciences, such as models for consensus formation, pedestrian dynamics and sociodynamics, as well as multispecies models in ecology \cite{valeria}. Our identification-control pipeline can also be applied to mean field models arising from rational strategies rather than interactions forces, e.g., mean field games with multiple solutions~\cite{Degond_al_2014}, or \color{black}mean field game models for systemic risk~\cite{carmona_2015}\color{black}. Moreover, the impact of the proposed methodology expands well beyond the realm of physical multiagent systems to problems in computational statistics, where there is interest in accelerating convergence to the target distribution for Fokker-Planck dynamics in applications related to sampling from probability measures in high dimensions~\cite{Kunisch_al_2018, Kunisch_al_2023, Lelievre_al_2013}.

\medskip

\paragraph{Acknowledgments:} GP is partially supported by an ERC-EPSRC Frontier Research Guarantee through Grant No. EP/X038645, ERC through Advanced Grant No. 247031. DK is partially supported by the EPSRC Standard Grant EP/T024429/1. SB has been supported by a Roth Scholarship from Imperial College London.

\section{Appendix}\label{sec:appendix}
\subsection{Fréchet differentiability of the cost functional  $\mathcal{J}$}\label{differentiability} \color{black}
This section examines the differentiability of the cost functional \eqref{ocp2}, which serves as the objective in our optimal control problem. The differentiability assumption is fundamental to the derivation of first-order optimality conditions, upon which our model predictive control algorithm is built.

By definition, $\mathcal{J}$ is Fréchet differentiable at $u\in \mathcal{U}:=L_2(0,T;L_\infty(\Omega))$ if there exists a bounded linear operator $D\mathcal{J}(u)[\cdot]$ such that 
\begin{equation}\label{frechet}
    \lim_{\|h\|_\mathcal{U}\to0}\frac{\big| \mathcal{J}(u+h) - \mathcal{J}(u) - D\mathcal{J}(u)[h]\big|}{\|h\|_\mathcal{U}} = 0,
\end{equation}
then $D\mathcal{J}(u)[h]$ is the Fréchet derivative of $\mathcal{J}$. 

The proof follows $3$ steps assessing the existence of a candidate operator matching the definition \eqref{frechet}, its linearity, and boundedness.
\begin{enumerate}
    \item \textbf{Candidate derivative -} We start by showing that the candidate $D\mathcal{J}(u)[h]$ we used to derive the optimality system \eqref{PDE_optimality} satisfies \eqref{frechet}. The expression for $D\mathcal{J}(u^*)[\cdot]$ was derived following \cite{PMP_PDE_troltzsch} and reads:
    $$D\mathcal{J}(u)[h]:=\int\limits_0^T\int\limits_\Omega (\rho(t,x)-\rho_\infty(x))\delta\rho(t,x) + \gamma\,u(t,x)\,h(t,x)\,dx\,dt\,,$$
    where $\delta\rho$ denotes the change in the state $\rho$ induced by the perturbation $u\mapsto u+h$.
    By expanding the perturbed cost functional $\mathcal{J}(u+h)$ we obtain:
     $$\begin{aligned}
       \mathcal{J}(u+h) &= \dfrac{1}{2}\int\limits_0^T\int\limits_\Omega (\rho - \rho_\infty + \delta\rho)^2 +\gamma(u+h)^2 dx\,dt\\
        &= \mathcal{J}(u) + \dfrac{1}{2} \int\limits_0^T\int\limits_\Omega 2\delta\rho(\rho-\rho_\infty) + \delta\rho^2+\gamma h^2 + 2\gamma u h\,dx\,dt\,.
    \end{aligned}$$
    Hence 
    $$\begin{aligned}
        \mathcal{J}(u+h) - \mathcal{J}(u) - D\mathcal{J}(u)[h] = 
        \dfrac{1}{2}\int\limits_0^T \|\delta\rho\|^2 + \gamma\|h\|^2\,dt\,,
    \end{aligned}$$
    where a standard energy estimate gives 
    \begin{equation}\label{drho}
        \|\delta\rho\|_{L_2(0,T;L_2(\Omega))} \leq C\|h\|_{\mathcal{U}}\,.
    \end{equation} 
    We conclude upon noticing that
    \begin{equation}\label{last_step}
        \frac{\big|\mathcal{J}(u+h) - \mathcal{J}(u) - D\mathcal{J}(u)[h]\big|}{\|h\|_\mathcal{U}}\;\leq\;\frac{\dfrac{(C+\gamma)}{2}\big\|h\big\|^2_\mathcal{U}}{\|h\|_\mathcal{U}} \xrightarrow[\|h\|_\mathcal{U}\mapsto0]{}0\,,
    \end{equation}
    hence the definition of Fréchet differential is satisfied.
    \item \textbf{Linearity -} The differential operator $D\mathcal{J}(u)[h]$ features both $h$ and the associated state perturbation $\delta\rho$. To show linearity in the argument $h$, we need to investigate the relationship between $h$ and $\delta\rho$. To this end, we derive the Fokker Planck equation for the state perturbation:
    \begin{equation}\label{linearizedFP}
        \partial_t \delta\rho = \beta^{-1} \Delta \delta\rho + \nabla \cdot (\delta\rho((\nabla W * \rho) + \nabla V + u)) + \nabla \cdot (\rho h)\,,
    \end{equation}
    which is linear in both $\delta\rho$ and $h$, so we can rewrite it as $\partial_t \delta\rho = \mathcal{L}(\delta\rho) + \nabla \cdot (\rho h)$, where $\mathcal{L}$ is a linear operator.    
    We introduce the solution operator $S$, so that $\delta\rho = S(h)$, which we will now show to be linear. Let $\delta\rho_1 = S(h_1)$ and $\delta\rho_2 = S(h_2)$, and consider $\delta\rho_{1,2} = a \delta\rho_1 + b \delta\rho_2$, for $a,b\in\R$. Substituting into the PDE \eqref{linearizedFP}, by linearity of $\mathcal{L}$ and $\nabla\cdot$ we have 
        \begin{align*}
\partial_t \delta\rho_{1,2} = \mathcal{L}(\rho_{1,2}) + \nabla \cdot (\rho(ah_1 + bh_2)),\,
\end{align*}
hence $\delta\rho_{1,2}$ satisfies the PDE for the control perturbation $a h_1 + b h_2$ and $S(a h_1 + b h_2) = a S(h_1) + b S(h_2)$. This implies that $\delta\rho$ depends linearly on $h$.
We can conclude that $D\mathcal{J}(u)[h]$ is indeed linear in $h$, as both $h$ and $\delta\rho$ appear linearly in the expression for $D\mathcal{J}(u)[h]$.

    \item \textbf{Boundedness -} We need to show that there exists a $ C>0$ such that 
    $$| D\mathcal{J}(u)[h] | \leq C\|h\|_\mathcal{U}.$$
\end{enumerate}
This reads
$$|D\mathcal{J}(u)[h] | \leq \bigg|\int\limits_{0}^T\int\limits_\Omega (\rho-\rho_\infty)\,\delta\rho\, dx \bigg| + \gamma\bigg| \int\limits_{0}^T\int\limits_\Omega u\,h\, dx \bigg|\,.$$
By Cauchy-Schwarz inequality and the estimate \eqref{drho}, we can bound the first term as 
\begin{align*}
    \bigg|\int\limits_{0}^T\int\limits_\Omega (\rho-\rho_\infty)\,\delta\rho\, dx \bigg|&\leq\|\rho-\rho_\infty\|_{L_2(0,T;L_2(\Omega))}\,\|\delta\rho\|_{L_2(0,T;L_2(\Omega))}\\&\leq\|\rho-\rho_\infty\|_{L_2(0,T;L_2(\Omega))}\,C\|h\|_{\mathcal{U}}\,\\[2ex] &\leq C_1 \|h\|_\mathcal{U}
\end{align*}
due to the regularity of $\rho$ and $\rho_\infty$. Similarly, the second term reads
\begin{equation*}
    \gamma\bigg| \int\limits_{0}^T\int\limits_\Omega u\,h\, dx \bigg| \leq \gamma \|u\|_{\mathcal{U}}\|h\|_{\mathcal{U}} \leq C_2 \|h\|_{\mathcal{U}}\,.
\end{equation*}
Combining the estimates for both terms, we obtain
$$|D\mathcal{J}(u)[h] |\leq (C_1+C_2)\|h\|_\mathcal{U}$$
which concludes the proof.

\subsection{Consistency of the discretized optimality condition}\label{consistency_appendix}
In this section we prove consistency between the two approaches to solving the PDE-constrained optimal control problem \eqref{ocp}. The crucial distinction lies at the stage at which discretization is applied to the solution process:
\begin{itemize}
    \item \textbf{optimize-then-discretize} - The optimality system is derived at the continuous, infinite-dimensional level, as in \eqref{PDE_optimality}. The functions $\rho$, $u$, and $\lambda$ are then discretized and the two point boundary value problem is rewritten in finite-dimensional form.
    \item \textbf{discretize-then-optimize} - Discretization is performed at the PDE level, resulting in a finite-dimensional optimal control problem. The optimality conditions for this approach are given by (\ref{state} - \ref{TPBVP}).
\end{itemize}
Consistency between these approaches can be shown by comparing the first order optimality conditions in \eqref{TPBVP} and \eqref{PDE_optimality}. 
The explicit expression of \eqref{TPBVP} reads
\begin{equation*}
     \gamma\, \mathbf{u}_j^\top M - \mathbf{p}^{\top} g_j(\mathbf{a}) = \gamma\int\limits {u}_{j}^L\,\psi_m\,d\Omega-\,\frac{1}{\int \psi_m \psi_m d\Omega} \int \rho^L \,\dfrac{\partial\lambda^L}{\partial{x_j}}\, \psi_m\,d\Omega 
\end{equation*}
for $j = 1,\dots, d$ and $m =1,\dots,L$ and $\lambda^L(t,x):=\sum_{n=1}^L \lambda_n(t)\psi_n(x)$. 
Noting that
\begin{itemize}
    \item since $u(t,\cdot)\in L_2(\Omega)\subset L_\infty(\Omega)$, and $\rho,\lambda\in \mathcal{W}(0,T)$ implies $\rho(t,\cdot)\in L_2(\Omega)$, $\nabla\lambda\in L_2(\Omega)$, as a direct result of Paserval's identity we have convergence of the truncated Fourier series $\rho^L,\,u^L,\,\frac{\partial\lambda^L}{\partial{x_j}}$;
    \item $\frac{1}{\int \psi_m \psi_m d\Omega}$ is constant and independent of the number of Fourier modes $L$, hence it represents a scaling that does not affect consistency;
\end{itemize}
we have 
$$ \gamma\, \mathbf{u}_j^\top M - \mathbf{p}^{\top} g_j(\mathbf{a}) \xrightarrow[L_2]{} \int (-\nabla\lambda\rho +\gamma u)\psi_m\, d\Omega$$
hence the discrete optimality condition converges in $L_2$ in variational sense to it's continuous counterpart, ensuring consistency of the discretization scheme.

\color{black}
\bibliographystyle{RS}
\bibliography{references}
\end{document}